\documentclass[12pt]{amsart}

\usepackage{latexsym}
\usepackage{amsfonts}
\usepackage{amssymb}
\usepackage{amsmath}
\usepackage{mathrsfs}
\usepackage{bbm}
\usepackage{dsfont}
\usepackage{hyperref}
\allowdisplaybreaks

\newcommand{\be}{\begin{equation}}
\newcommand{\ee}{\end{equation}}
\newcommand{\bea}{\begin{eqnarray}}
\newcommand{\eea}{\end{eqnarray}}
\newcommand{\bean}{\begin{eqnarray*}}
\newcommand{\eean}{\end{eqnarray*}}
\newcommand{\brray}{\begin{array}}
\newcommand{\erray}{\end{array}}

\newtheorem{dfn}{Definition}[section]
\newtheorem{thm}[dfn]{Theorem}
\newtheorem{lmma}[dfn]{Lemma}
\newtheorem{ppsn}[dfn]{Proposition}
\newtheorem{crlre}[dfn]{Corollary}
\newtheorem{xmpl}[dfn]{Example}
\newtheorem{rmrk}[dfn]{Remark}

\newcommand{\bdfn}{\begin{dfn}\rm}
\newcommand{\bthm}{\begin{thm}}
\newcommand{\blmma}{\begin{lmma}}
\newcommand{\bppsn}{\begin{ppsn}}
\newcommand{\bcrlre}{\begin{crlre}}
\newcommand{\bxmpl}{\begin{xmpl}}
\newcommand{\brmrk}{\begin{rmrk}\rm}

\newcommand{\edfn}{\end{dfn}}
\newcommand{\ethm}{\end{thm}}
\newcommand{\elmma}{\end{lmma}}
\newcommand{\eppsn}{\end{ppsn}}
\newcommand{\ecrlre}{\end{crlre}}
\newcommand{\exmpl}{\end{xmpl}}
\newcommand{\ermrk}{\end{rmrk}}

\newcommand{\clh}{\mathcal{H}}

\newcommand{\clk}{\mathcal{K}}

\author{Piyasa Sarkar and S. Sundar}
\title {Examples of  multiparameter CCR flows with non-trivial index}

\begin{document}
\maketitle

\begin{abstract}
In this paper, we construct uncountably many examples of multiparameter CCR flows, which are not pullbacks of $1$-parameter CCR flows, with index one. Moreover, the constructed CCR flows are type I in the sense that the associated product system is the smallest subsystem containing  its units. 
\end{abstract}

\noindent {\bf AMS Classification No. :} {Primary 46L55; Secondary 46L99.}  \\
{\textbf{Keywords :}} $E_0$-semigroups, CCR flows, Index. 

\section{Introduction}
The remarkable efforts of R.T. Powers, Arveson, Tsirelson and many others have firmly established  the study of one parameter $E_0$-semigroups  on a type I factor as a fruitful area of research during the last thirty years.  Arveson's monograph \cite{Arveson} forms the definitive reference for the subject. 
 Attempts to extend the theory of $E_0$-semigroups/CP-semigroups to the multiparameter context have been made in the recent years.    A few  recent  papers that explore these issues are  \cite{Shalit_2008}, \cite{Shalit_2011}, \cite{Shalit_Solel}, \cite{Shalit} ,\cite{Skeide_Shalit} \cite{Anbu}, \cite{Anbu_Sundar}, and \cite{Murugan_Sundar_continuous}. This article is one such attempt in this direction.

Consider a pointed, spanning, closed convex cone $P$ in $\mathbb{R}^{d}$.  An $E_0$-semigroup  over $P$  is a semigroup, indexed by $P$, of unital normal $*$-endomorphisms of $B(\clh)$ which is continuous in an appropriate sense. 
The first basic example in the theory is that of a CCR flow.   We can associate to each isometric representation $V$ of $P$ on a Hilbert space $\clh$, an $E_0$-semigroup, denoted $\alpha^{V}$ and called the CCR flow associated to $V$, on the algebra of bounded operators of the symmetric Fock space $\Gamma(\clh)$. 

In the one parameter case, i.e. when $P=[0,\infty)$,   it is well known from the seminal work of Arveson (\cite{Arv_Fock}) that an $E_0$-semigroup is a CCR flow if and only if it has units and the units generate the associated product system. 
It follows from the   results of Tsirelson (\cite{Tsi}, \cite{Tsirelson2003}), Powers (\cite{Powers_TypeIII}, \cite{Powers}), Liebscher (\cite{Liebscher}) and Izumi and Srinivasan (\cite{Vasanth_Izumi}, \cite{Vas_Izumi}),
that  other types of $E_0$-semigroups exist, i.e.   there exist  
 uncountably many one parameter $E_0$-semigroups which have units but whose units do not generate the associated product system, and
  uncountably many one parameter $E_0$-semigroups which do not have a unit. 

Unlike in the one parameter situation, the first examples of multiparameter CCR flows considered in \cite{Anbu} and   \cite{Anbu_Sundar} admit only one unit, up to a character and  it is the vaccum unit. This is because, up to characters, units of CCR flows are in bijective correspondence with additive cocycles of the associated isometric representations and the isometric representations considered in \cite{Anbu_Sundar} fail to have a  non-zero additive cocycle (see Prop. 2.4 of \cite{Anbu_Sundar}). This is a bit disconcerting as the characteristic feature of one parameter CCR flows is that they have units in abundance. 

Given this anomaly, the role of units in the one parameter theory and the  literature alluded to above, it  is natural and is of intrinsic interest  to ask whether, in the higher dimensional case, there exist CCR flows which have a unit other than the vaccum unit.  
Multiparameter CCR flows with more than one unit can easily be constructed by pulling back a  one parameter CCR flow by a homomorphism $\phi:P \to [0,\infty)$. However, this is tautological  and we do not consider such examples as multiparametric in nature.   A more pertinent question is the following.

\vspace{2.5 mm}
\noindent 
 \textbf{Question 1:} In the higher dimensional case, are there  CCR flows, which are not obtained as pullbacks of one parameter CCR flows, with more than one unit and  which is also ``type I" in some sense ? 
\vspace{2.5 mm}

 We construct uncountably many such CCR flows with Arveson's index one.   We also prove that the associated product system does not have a  smaller subsystem containing all its units. In this sense, the constructed CCR flows are type I. 

Next, we explain the strategy used  to construct such examples. A $P$-space, say $A$, in $\mathbb{R}^{d}$ is a proper non-empty closed  subset which is $P$-invariant, i.e. $A+x \subset A$ for $x \in P$. If $\clk$ is a Hilbert space of dimension $k$, the translation action of $P$ on $A$, implements an isometric representation of $P$ on $L^{2}(A,\clk)$. The corresponding CCR flows $\alpha^{(A,k)}$ were investigated in greater detail in \cite{Anbu_Sundar}. 

Here, rather than looking at $P$-spaces in $\mathbb{R}^{d}$, we look at $P$-spaces $A$ in the quotient group $G:=\mathbb{R}^{d}/H$ where $H$ is a closed subgroup of $\mathbb{R}^{d}$. Let $\widetilde{V}^{(A,k)}$ be the isometric representation of $P$ on $L^{2}(A,\clk)$ corresponding to the translation action of $P$ on $A$. Denote the associated CCR flow by $\widetilde{\alpha}^{(A,k)}$. We show that the CCR flow $\widetilde{\alpha}^{(A,k)}$  has more than one unit if and only if  $H$ has codimension one.\footnote{Let $V$ be the maximal vector subspace of $H$. We say that $H$ has codimension one  if the discrete subgroup $H/V$ in $\mathbb{R}^{d}/V$ has rank $d-1-dim(V)$.} 

On closer inspection, we notice that  $G$ is  an abelian Lie group and $Q:=\overline{\pi(P)}$ is an abelian Lie semigroup. Here $\pi$ is the quotient map. Moreover, $A$ is a $Q$-space. Denote the CCR flow, indexed by $Q$, associated to the shift action of $Q$ on $L^{2}(A,\clk)$ by $\alpha^{(A,k)}$. Then, the $E_0$-semigroup
$\widetilde{\alpha}^{(A,k)}$  is  the pull back of the $E_0$-semigroup $\alpha^{(A,k)}$ via the map $\pi:P \to Q$. Thus, ultimately the analysis of $\widetilde{\alpha}^{(A,k)}$ depends on the analysis of  the CCR flows $\alpha^{(A,k)}$, indexed by $Q$.
After having come this far, it is only appropriate to analyse the CCR flows $\alpha^{(A,k)}$ indexed by a Lie semigroup, say $P \subset G$, where $A$ is a $P$-space in $G$.

The second motivation for us to consider arbitrary non-commutative Lie semigroups comes from the main result (Theorem 1.2) of \cite{Anbu_Sundar} which asserts that, in the case of a cone, the CCR flow $\alpha^{(A,k)}$ remembers both $A$, up to a translate,  and the multiplicity $k$. The proof makes  heavy weather of the fact that the group law involved is abelian (Prop. 4.4, Lemma 4.6 of \cite{Anbu_Sundar}) and the proof works only in the commutative situation. 
 A cone being a prototypical example of a Lie semigroup, it is of interest to know whether the same result is valid in the non-Euclidean and in the non-commutative situtation. 
The second question that we investigate in this paper is the following. 

\vspace{2.5 mm}
\noindent
 \textbf{Question 2:} Let $P$ be a Lie semigroup, not necessarily commutative. Does the CCR flow $\alpha^{(A,k)}$ remember the set $A$, up to a translate, and the multiplicity $k$ ?
\vspace{2.5 mm}

The affirmative answer to the above question for the case of a cone was first given in  \cite{Anbu_Sundar}. The method used in \cite{Anbu_Sundar} relies  on computing the gauge groups of the associated CCR flows and the groupoid machinery developed in \cite{Sundar_Ore}.  Later in \cite{injectivity}, a better conceptual explanation clarifying the  exact role played by groupoids was obtained (see Theorem 5.2 of \cite{injectivity} and the discussion following it).  This, in turn, depends on establishing the result that, for a pure isometric representation $V$,  the  representation $V$ can be recovered from the cocycle conjugacy class of the CCR flow $\alpha^{V}$.

However, the proof given in \cite{injectivity}  is very long and not quite direct.  A shorter proof was found by R. Srinivasan in \cite{Vasanth}. In this paper, we give yet another direct proof of this result in the setting of Lie semigroups. We then apply this to settle Question 2. We show that, in general, for an arbitrary Lie semigroup $P$, the CCR flow $\alpha^{(A,k)}$ remembers the set $A$ but not necessarily the multiplicity $k$.  This is in stark contrast to the case of a cone.

The organisation of this paper is next described. 

 After this introductory section, in Section 2, we collect the preliminaries on Lie semigroups and $E_0$-semigroups that are required to read this paper. Imitating Arveson, when $P$ is a closed convex cone, we define the notion of index, a numerical invariant, which measures the ``number of units" of a spatial $E_0$-semigroup. We show that for a CCR flow index coincides with the dimension of the space of additive cocycles. 

In Section 3, we calculate the units of the CCR flow associated to a $P$-space of multiplicity $1$. We show in particular that if the enveloping group is unimodular, then the CCR flow associated to a $P$-space $A$ has more than one unit if and only if the boundary of $A$ is compact. In Section 4, we derive a necessary condition for the boundary of a $P$-space to be  compact which we prove is also sufficient in the abelian case. In particular, we show that if $P$ is a generating Lie semigroup of an abelian Lie group of the form $\mathbb{R}^{d} \times \mathbb{T}^{s}$ and $A$ is a $P$-space, then the boundary of $A$ is compact if and only if $d=1$. 

In Section 5, we provide another direct proof of the fact that  for two pure isometric representations $V$ and $W$, the corresponding CCR flows $\alpha^{V}$ and $\alpha^{W}$ are cocycle conjugate if and only if $V$ and $W$ are unitarily equivalent.  We prove this in the setting of Lie semigroups. We apply this to study the CCR flows $\alpha^{(A,k)}$. We avoid groupoids completely in this paper and instead use a simple dilation trick. We show that $\alpha^{(A,k)}$ remembers $A$, up to a translate, but not necessarily the multiplicity $k$. 
In Section 6, for a pointed and a spanning closed convex cone in $\mathbb{R}^{d}$, where $d \geq 2$, we exhibit uncountably many type I examples of CCR flows with index one that are not pullbacks of one parameter CCR flows.

\section{Preliminaries}
In this section, we collect the preliminaries required to read the rest of this paper. Let $G$ be a  connected Lie group with Lie algebra $\mathfrak{g}$. Denote the exponential map by $\exp$. Let $P$ be a closed subsemigroup of $G$ containing the identity element $e$. The Lie wedge of $P$, denoted $L(P)$, is defined as follows.
\[
L(P):=\{X \in \mathfrak{g}: exp(tX) \in P \textrm{~for all $t \geq 0$}\}.\]
The semigroup $P$ is called \emph{a Lie semigroup} if the semigroup generated by $\exp(L(P))$ is dense in $P$. 
Throughout this paper, we assume that the Lie semigroups that we consider have dense interior. For a Lie semigroup $P$, we denote its interior by $\Omega$. The condition that $\Omega$ is dense in $P$ is equivalent to the condition that the Lie algebra generated by $L(P)$ is  $\mathfrak{g}$ (see Corollary 5.12 of \cite{Neeb_Hilgert}). Throughout this paper, we exclude the case $P=G$. 

For the rest of  this section, we assume that $P$ is a Lie semigroup in $G$ with dense interior $\Omega$. We also assume that $P$ is Ore in $G$, i.e. $PP^{-1}=P^{-1}P=G$. For $x,y \in G$, we say $x \leq y$ ($x<y$) if $x^{-1}y \in P$ ($x^{-1}y \in \Omega$). Note that $\leq$ is a preorder and is a partial order if and only if $P \cap P^{-1}=\{e\}$. The following are some examples that fit our assumptions. 

\begin{enumerate}
\item[(1)] Let $G=\mathbb{R}^{d}$ and let $P$ be a closed convex cone in $\mathbb{R}^{d}$. 
\item[(2)] \textbf{The $ax+b$-semigroup:} Let \[G:=\Big\{\begin{bmatrix}
               a & b \\
               0 & 1
               \end{bmatrix}: a>0, b \in \mathbb{R}\Big\}\] and let \[P:=\Big\{\begin{bmatrix}
                                                                                                         a & b \\
                                                                                                         0 & 1
                                                                                                         \end{bmatrix}: a \geq 1, b \geq 0 \Big\}.\]
\item[(3)] \textbf{The Heisenberg semigroup:}         Let $G:=H_{2n+1}=\mathbb{R}^{n} \times \mathbb{R}^{n} \times \mathbb{R}$. 
 The group law on $H_{2n+1}$ is defined by
 \[
 (x_1,y_1,z_1)(x_2,y_2,z_2):=(x_1+x_2,y_1+y_2,z_1+z_2+\langle x_1|y_2 \rangle).\] 
  Here $\langle~|~\rangle$ stands for the usual inner product.       Let \[P:=H_{2n+1}^{+}:=\{(x,y,z) \in H_{2n+1}: x \in \mathbb{R}_{+}^{n}, y \in \mathbb{R}_{+}^{n}, z \in \mathbb{R}_{+}\}.\]                                                                                        
\item[(4)] Let $G$ be the group of  upper triangular real matrices with $1$'s in the diagonal and let $P$ be the subsemigroup of $G$ whose entries are non-negative. 
\item[(5)] \textbf{The contraction semigroup:} Let $V$ be a finite dimensional inner product space. Set $G:=GL(V)$ and let $P:=\{T \in G: ||T|| \leq 1\}$. 
\end{enumerate}
\begin{rmrk}
\label{interior}
An important fact that we will repeatedly use and is worth stating is the following (see Corollary 3.11 of \cite{Neeb_Hilgert}). 
 Let $a \in \Omega$ be given. Then there exist a natural number $n \geq 1$ and $X_1,X_2,\cdots,X_n \in L(P)$ such that 
$
a=\exp(X_1)\exp(X_2)\cdots \exp(X_n)$. 
\end{rmrk}

Next, we recall the basic definitions concerning the theory of $E_0$-semigroups.
\begin{dfn}
Let $\clh$ be a separable infinite dimensional Hilbert space. An $E_0$-semigroup over $P$ on $B(\clh)$ is a family $\alpha:=\{\alpha_x\}_{x \in P}$ of unital normal $*$-endomorphisms of $B(\clh)$ such that 
\begin{enumerate}
\item[(1)] for $x,y \in P$, $\alpha_x \circ \alpha_y=\alpha_{xy}$, and
\item[(2)] for $A \in B(\clh)$ and $\xi,\eta \in \clh$, the map $P \ni x \to \langle \alpha_x(A)\xi|\eta \rangle \in \mathbb{C}$ is continuous. 
\end{enumerate}
\end{dfn}

Let $\alpha:=\{\alpha_x\}_{x \in P}$ and $\beta:=\{\beta_x\}_{x \in P}$ be two $E_0$-semigroups over $P$ on $B(\clh)$. We say that $\beta$ is a \emph{cocycle perturbation} of $\alpha$ if there exists a strongly continuous family of unitaries $\{U_x\}_{x \in P}$ in $B(\clh)$ such that 
\begin{enumerate}
\item[(a)] for $x \in P$ and $A \in B(\clh)$, $\beta_x(A)=U_x\alpha_x(A)U_x^{*}$, and
\item[(b)] for $x,y \in P$, $U_{xy}=U_x\alpha_x(U_y)$. 
\end{enumerate}
Let $\alpha:=\{\alpha_x\}_{x \in P}$ and $\beta:=\{\beta_x\}_{x \in P}$ be $E_0$-semigroups over $P$ on $B(\clh)$ and $B(\clk)$ respectively. We say that $\alpha$ and $\beta$ are \emph{cocycle conjugate} if for a  unitary $U:\clh \to \clk$, $\beta$ is a cocycle perturbation of $\{Ad(U) \circ \alpha_x \circ Ad(U^{*})\}_{x \in P}$. 
 When $P$ is clear from the context, we simply call an $E_0$-semigroup over $P$ an $E_0$-semigroup. 
 
The $E_0$-semigroups that we study in  this paper are CCR flows. Let $V:=\{V_x\}_{x \in P}$ be a strongly continuous semigroup of isometries, also called \emph{an isometric representation of $P$}, on a Hilbert space $\clh$. Denote the symmetric Fock space of $\clh$ by $\Gamma(\clh)$. Then there exists a unique $E_0$-semigroup, denoted $\alpha^{V}:=\{\alpha_x\}_{x \in P}$, on $B(\Gamma(\clh))$ such that for $x \in P$ and $\xi \in \clh$, 
\[
\alpha_x(W(\xi))=W(V_x\xi).\]
Here $\{W(\xi):\xi \in \clh\}$ denotes the collection of Weyl operators on $\Gamma(\clh)$. Recall that the action of the Weyl operators on the set of exponential vectors $\{e(\eta): \eta \in \clh\}$ is given by the  equation 
\[
W(\xi)e(\eta)=e^{-\frac{||\xi||^{2}}{2}-\langle \eta|\xi \rangle}e(\xi+\eta).\]
The $E_0$-semigroup $\alpha^{V}$ is called \emph{the CCR flow associated to the isometric representation $V$}.

The class of isometric representations that we will be interested in this paper is described below. Let $A \subset G$ be a non-empty proper closed subset. We say that $A$ is a $P$-space if $PA \subset A$, i.e $ax \in A$ for $a \in P$ and $x \in A$. Let $\clk$ be a Hilbert space of dimension $k$.  Consider the Hilbert space $\clh:=L^{2}(A,\clk)$. For $x \in P$, define the isometry $V_{x}$ on $L^{2}(A,\clk)$ by the following formula.  \begin{equation}
\label{isometries}
V_{x}(f)(y):=\begin{cases}
 f(x^{-1}y)  & \mbox{ if
} x^{-1}y \in A,\cr
   &\cr
    0 &  \mbox{ if } x^{-1}y \notin A
         \end{cases}
\end{equation}
for $f \in L^{2}(A,\clk)$. Then $V:=\{V_{x}\}_{x \in P}$ is a strongly continuous isometric representation of $P$ on $L^{2}(A,\clk)$. Also, the representation $V$ is pure, i.e. $\displaystyle \bigcap_{x \in P}Ran(V_x)=\{0\}$ (see Example 11.8 of \cite{Sundar_Notes}).  We call $V$ \emph{the isometric representation associated to the $P$-space $A$ of multiplicty $k$} and denote it by $V^{(A,k)}$.  We denote the CCR flow associated to the isometric representation $V^{(A,k)}$ by $\alpha^{(A,k)}$. 

Let $\alpha:=\{\alpha_x\}_{x \in P}$ be an $E_0$-semigroup on $B(\clh)$. Suppose $u:=\{u_x\}_{x \in P}$ is a strongly continuous family of bounded operators on $\clh$. We call $u$  \emph{a unit of $\alpha$} if 
\begin{enumerate}
\item[(1)] for $x \in P$ and $A \in B(\clh)$, $\alpha_x(A)u_x=u_xA$, 
\item[(2)] for $x \in P$, $u_x \neq 0$, and
\item[(3)] for $x,y \in P$, $u_{xy}=u_xu_y$. 
\end{enumerate}
Denote the set of units of $\alpha$ by $\mathcal{U}_{\alpha}$. The $E_0$-semigroup $\alpha$ is said to be  \emph{spatial} if $\mathcal{U}_{\alpha} \neq \emptyset$. For a CCR flow $\alpha$, $\mathcal{U}_{\alpha}$ has a neat description which we explain below. 

Let $V$ be an isometric representation of $P$ on a Hilbert space $\clh$. A continuous map $\xi:P \to \clh$ is said to be \emph{an additive cocycle of $V$} if 
\begin{enumerate}
\item[(1)] for $x \in P$, $\xi_x \in Ker(V_x^{*})$, and
\item[(2)] for $x,y \in P$, $\xi_{xy}=\xi_x+V_x\xi_y$. 
\end{enumerate}
Denote the vector space of additive cocycles of $V$ by $\mathcal{A}(V)$. 

Fix an isometric representation $V$ of $P$ on a Hilbert space $\clh$. Let $\alpha$ be the CCR flow associated to $V$. 
Let $x \in P$ and $\xi \in Ker(V_x^{*})$ be given. Define a bounded operator, denoted $T^{(x)}_{e(\xi)}$, on $\Gamma(\clh)$ by the following equation. 
\[
T^{(x)}_{e(\xi)}(e(\eta))=e(\xi+V_x\eta).\]
 Denote the set of continuous group homomorphisms from $G$ to the multiplicative group $\mathbb{C}^{\times}:=\mathbb{C} \backslash \{0\}$ by $Hom(G,\mathbb{C}^{\times})$.  Let  $\chi \in Hom(G,\mathbb{C}^{\times})$ and $\xi \in \mathcal{A}(V)$ be given. For $x \in P$, let \[
 u^{(\chi,\xi)}_x:=\chi(x)T^{(x)}_{e(\xi_x)}.\]
 It is easy to verify that  $u^{(\chi,\xi)}:=\{u^{(\chi,\xi)}_{x}\}_{x \in P}$ is a unit of $\alpha$. 
 
 With the foregoing notation, we have the following proposition. Since the proof is similar to Theorem 5.10 of \cite{Anbu}, we merely give a sketch of the proof. 
 
 \begin{ppsn}
 \label{units}
 The map 
 \[
 Hom(G,\mathbb{C}^{\times}) \times \mathcal{A}(V) \ni (\chi,\xi) \to u^{(\chi,\xi)} \in \mathcal{U}_{\alpha}\]
 is a bijection. 
  \end{ppsn}
 \textit{Proof.} The only point that requires explanation is the surjectivity of the map. Let $u:=\{u_x\}_{x \in P}$ be a unit of $\alpha$. 
 Set 
 \[
 Q:=\{x \in P: \textrm{there exist $\lambda_x \in \mathbb{C}^{\times}$ and $\xi_x \in Ker(V_x^{*})$ such that $u_x=\lambda_xT^{(x)}_{e(\xi_x)}$}\}.\]
Note that $Q$ is a closed subsemigroup of $P$. 

Fix $X \in L(P)$. For $t \geq 0$, let $\alpha_t:=\alpha_{\exp(tX)}$. Then $\{\alpha_t\}_{t \geq 0}$ is the one parameter CCR flow associated to the isometric representation $\{V_{\exp(tX)}\}_{t \geq 0}$. Making use of Wold decomposition and Theorem 2.6.4 of \cite{Arveson}, it is routine to see that $\exp(tX) \in Q$ for every $t \geq 0$. Hence, $Q$ contains the semigroup generated by $\exp(L(P))$. But the semigroup generated by $\exp(L(P))$ is dense in $P$. Therefore $Q=P$. 

Thus, for every $x \in P$, there exists $\lambda_x \in \mathbb{C}^{\times}$ and $\xi_x \in Ker(V_x^{*})$ such that \[u_x=\lambda_xT^{(x)}_{e(\xi_x)}.\] The rest of the proof is similar to Theorem 5.10 of \cite{Anbu}. \hfill $\Box$

 \textbf{Arveson's index:} 
 If $P$ is a closed convex cone then, just like in the $1$-dimensional case, a numerical invariant, called \emph{the index}, can be defined which measures the relative abundance of units of a spatial $E_0$-semigroup. For the rest of this section, assume that $P$ is a closed convex cone in $\mathbb{R}^{d}$ with dense interior $\Omega$. 
  Let $\alpha:=\{\alpha_x\}_{x \in P}$ be an $E_0$-semigroup over $P$ on $B(\clh)$. Assume that $\alpha$ is spatial. Denote the collection of units by $\mathcal{U}_{\alpha}$. 
  
  Fix $a \in \Omega$. It is easy to show that for $u,v \in \mathcal{U}_{\alpha}$, the map $(0,\infty) \ni t \to v_{ta}^{*}u_{ta} \in \mathbb{C}\backslash\{0\}$ is multiplicative. Thus, for $u,v \in \mathcal{U}_{\alpha}$, there exists a unique complex number $c_{a}(u,v)$ such that for $t>0$,
\[
v_{ta}^{*}u_{ta}=e^{tc_{a}(u,v)}.\]
The map $c_{a}:\mathcal{U}_{\alpha} \times \mathcal{U}_{\alpha} \to \mathbb{C}$ is conditionally positive definite. Let $\clh(\mathcal{U}_{\alpha})$ be the Hilbert space obtained by applying the usual GNS construction. 

Let us recall the construction of $\clh(\mathcal{U}_{\alpha})$. Denote by $\mathbb{C}_{0}(\mathcal{U}_{\alpha})$ the space of finitely supported functions on $\mathcal{U}_{\alpha}$ whose sum is zero. For $f,g \in \mathbb{C}_{0}(\mathcal{U}_{\alpha})$, define
\[
\langle f|g \rangle_a:=\sum_{u,v \in \mathcal{U}_{\alpha}}f(u)\overline{g(v)}c_{a}(u,v).\]
Then $\langle~|~\rangle_a$ is a semi-definite inner product on $\mathbb{C}_{0}(\mathcal{U}_{\alpha})$. The Hilbert space $\clh(\mathcal{U}_{\alpha})$ is obtained by completing the quotient of  $\mathbb{C}_{0}(\mathcal{U}_{\alpha})$ by the space of null vectors. The Hilbert space $\mathcal{H}(\mathcal{U}_{\alpha})$ is separable. For units $u,v \in \mathcal{U}_{\alpha}$, we denote the class representing $\delta_u-\delta_v$ by $[u]-[v]$. Here, $\delta_u$ stands for the characteristic function at $u$.

\begin{ppsn}
Keep the foregoing notation. The dimension of the Hilbert space $\mathcal{H}(\mathcal{U}_{\alpha})$ is independent of the chosen interior point $a$. 
\end{ppsn}
\textit{Proof.} Observe that for units $u,v \in \mathcal{U}_{\alpha}$, the map $P \ni x\to u_x^{*}v_x \in \mathbb{C}^{\times}$ is multiplicative. Let $a,b \in \Omega$ be given. Denote the seminorms corresponding to the semi-definite inner products $\langle~|~\rangle_a$ and $\langle~|~\rangle_b$ by $||~||_a$ and $||~||_b$ respectively.

Let $b \in \Omega$. By the Archimedean property, there exists $n \in \mathbb{N}$ such that $b<na$, i.e. $na-b \in \Omega$. So, $na=b+c$, for some $c \in \Omega$. Let $u,v \in \mathcal{U}_{\alpha}$ be given. Then, 
\[ e^{tc_{na}(u,v)}=e^{ntc_{a}(u,v)}, \quad e^{tc_{b+c}(u,v)}=e^{t(c_{b}(u,v)+c_{c}(u,v))} \]
 which  gives us the relation $nc_{a}(u,v)=c_{b}(u,v)+c_{c}(u,v)$. Similarly, there exists $m \in \mathbb{N}$ such that $mb=a+d$, for some $d \in \Omega$, and $mc_{b}(u,v)=c_{a}(u,v)+c_{d}(u,v)$. Combining the two, we have for $f \in C_0(\mathcal{U}_\alpha)$, $n||f||_{a}^{2}=||f||_{b}^{2}+||f||_{c}^{2}$ and $m||f||_{b}^{2}=||f||_{a}^{2}+||f||_{d}^{2}$. Therefore, for $f \in C_0(\mathcal{U}_\alpha)$, 
 \[ ||f||_{a} \leq \sqrt{m}||f||_{b} \quad \textrm{and} \quad ||f||_{b} \leq \sqrt{n} ||f||_{a}. \]
 Thus, both the seminorms $||~||_a$ and $||~||_b$ are equivalent. Consequently, the dimension of $\mathcal{H}(\mathcal{U}_{\alpha})$ is independent of the chosen interior point. \hfill $\Box$

\begin{dfn}
Let $\alpha$ be a spatial $E_0$-semigroup over $P$.    The index of $\alpha$, denoted $Ind(\alpha)$, is defined to be the dimension of the Hilbert space $\mathcal{H}(\mathcal{U}_\alpha)$.  
  \end{dfn}
 
 For CCR flows, computing the index is equivalent to computing the dimension of the space of additive cocycles. Let $V$ be an isometric representation of $P$ on a Hilbert space $\clh$. Denote the CCR flow associated to $V$ by $\alpha$. The space of additive cocycles of $V$ will be denoted by $\mathcal{A}(V)$. Fix $a \in \Omega$. For $\xi,\eta \in \mathcal{A}(V)$, let 
 \[\langle \xi|\eta \rangle _a:=\langle \xi_a|\eta_a \rangle.\]
 
 \begin{lmma}
 The sesquilinear form $\langle~|~\rangle_a$ is an inner product on $\mathcal{A}(V)$ and $\mathcal{A}(V)$ is a Hilbert space with respect to $\langle~|~\rangle_a$. 
  \end{lmma}
  \textit{Proof.} Note that for $\xi \in \mathcal{A}(V)$, the map $P \ni x \to \langle \xi_x|\xi_x \rangle \in [0,\infty)$ is additive. Let $\xi \in \mathcal{A}(V)$ be given.  Let $x \in P$ and let $n$ be a natural number such that $x \leq na$. Then, $na=x+y$ for some $y \in P$. 
    Note that $\xi_{na}=\xi_{x}+V_{x}\xi_{y}$, and consequently
  $ n ||\xi_{a}||^{2}=\langle \xi_{na}|\xi_{na}\rangle=\langle \xi_x| \xi_x\rangle+\langle \xi_y|\xi_y \rangle$. 
  This implies \begin{equation}\label{inequality monotone} ||\xi_x||^2 \leq n ||\xi_a||^2.\end{equation} 
  
  By the Archimedean property, given any $x \in P$, we have an $n \in \mathbb{N}$ such that $na=x+y$, for some $y \in P$.  By Eq. \ref{inequality monotone}, if for some $\xi$, $\langle \xi | \xi \rangle _a= \langle \xi_a| \xi_a \rangle =0$, then $||\xi_x||=0$, for all $x \in P$, and $\xi=0$ in $\mathcal{A}(V)$. This proves that the sesquilinear form is indeed an inner product.
  
  Given $b \in \Omega$,  there exist natural numbers $k$ and $m$  such that $ka-b \in P$ and $mb-a \in P$. By Eq. \ref{inequality monotone},  
  $\frac{1}{k}||\xi_b||^2 \leq ||\xi_a||^2 \leq m ||\xi_b||^2 $.
  Thus, the norms induced by any two interior points $a$ and $b$ are equivalent.
  
  To prove that $\mathcal{A}(V)$ is a Hilbert space, let $\left\lbrace \xi^{n} \right\rbrace $ be a Cauchy sequence in $\mathcal{A}(V)$. Let $x \in P$ be given. Choose a natural number $k$ such that $ka>x$. Then, by Eq. \ref{inequality monotone}, for $n,m \geq 1$, 
  \[
  ||\xi^n_x-\xi^m_x|| \leq k ||\xi^n_a-\xi^m_a||.\]
     Hence, for all $x\in P$, $\left\lbrace \xi^{n}_{x}\right\rbrace $ is Cauchy in $Ker(V_x^{*})$ and hence converges to, say,  $\eta_{x}$ in $Ker(V_{x}^{*})$. 
  
  Fix $ x \in P$, $0< \delta < 1$ and $a \in \Omega$. By the Archimedean Property, given $y \in B(x, \delta) \cap P$, there exists $k \in \mathbb{N}$ such that $y < ka$. Now, $P=\displaystyle \bigcup_{n \geq 1} \left\lbrace y \in P | y <na\right\rbrace $ and $B(x,\delta) \cap P$ has compact closure. So, there exists $k \in \mathbb{N}$ such that $B(x, \delta) \cap P \subset \left\lbrace y \in P| y< ka\right\rbrace $. Since $y <ka$, for all $y \in B(x, \delta) \cap P$, we have by Eq. \ref{inequality monotone} that for all $n,m$, 
  \[
  ||\xi^n_y-\xi_y^m|| \leq k ||\xi^n_a-\xi^m_a||.\]
  This shows that the sequence $\left\lbrace \xi^n \right\rbrace $ is locally uniformly Cauchy.  Consequently, the map $\eta = \left\lbrace \eta_x \right\rbrace_{x \in P} $ is continuous. Also, $\eta_{x+y}$ being the limit of $\xi^{n}_{x+y}=\xi^{n}_{x}+V_{x} \xi^{n}_{y}$, we have for $x,y \in P$, $\eta_{x+y}=\eta_{x}+V_{x}\eta_{y}$. Thus, $\eta$ belongs to $\mathcal{A}(V)$. This completes the proof. \hfill $\Box$

  In the next proposition, we identify $Hom(\mathbb{R}^{d},\mathbb{C}^{\times})$ with $\mathbb{C}^{d}$. The identification is via the map 
  \[
  \mathbb{C}^{d} \ni \lambda \to (\mathbb{R}^{d} \ni x \to e^{\langle \lambda|x \rangle} \in \mathbb{C}^{\times}) \in Hom(\mathbb{R}^{d},\mathbb{C}^{\times}).\]
  \begin{ppsn}
  The map \[\mathcal{A}(V) \ni \xi \to [u^{(0,\xi)}]-[u^{(0,0)}] \in \mathcal{H}(\mathcal{U}_\alpha)\] is a unitary. Consequently, $Ind(\alpha)=\dim \mathcal{A}(V)$. 
    \end{ppsn}
  \textit{Proof.} Fix $a \in \Omega$. For two additive cocycles $\xi, \eta \in \mathcal{A}(V)$ and for any $\lambda , \mu \in \mathbb{C}^d$, note that \[c_{a}(u^{(\lambda,\xi)},u^{(\mu,\eta)}) = \langle \lambda|a \rangle + \overline{\langle \mu|a \rangle} + \langle \xi|\eta \rangle_a.\] Calculate as follows to observe that 
  \begin{align*}
  & \Big \langle [u^{(0,\xi)}]-[u^{(0,0)}] \Big| [u^{(0,\eta)}]-[u^{(0,0)}] \Big \rangle  \\ 
  &= c_{a}(u^{(0,\xi)},u^{(0,\eta)}) - c_{a}(u^{(0,\xi)},u^{(0,0)}) - c_{a}(u^{(0,0)},u^{(0,\eta)}) + c_{a}(u^{(0,0)},u^{(0,0)}) \\ 
  &= \langle \xi| \eta \rangle_a - \langle \xi |0 \rangle_a - \langle 0 | \eta \rangle_a + \langle 0|0 \rangle_a \\ &= \langle \xi | \eta \rangle_a.
    \end{align*} 
   Thus, the prescribed map preserves the inner product. A routine calculation shows that \[ \langle [u^{(\lambda,\xi)}]-[u^{(0,\xi)}] | [u^{(\lambda,\xi)}]-[u^{(0,\xi)}] \rangle =0.\] Thus, $ [u^{(\lambda,\xi)}]-[u^{(0,0)}]=[u^{(0,\xi)}]-[u^{(0,0)}]$. Also, the set $\left\lbrace [u]-[v]| u, v \in \mathcal{U}_{\alpha} \right\rbrace $ is total in $\mathcal{H}(\mathcal{U}_\alpha)$ and \[[u^{(\lambda,\xi)}]-[u^{(\mu,\eta)}]=([u^{(0,\xi)}]-[u^{(0,0)}])-([u^{(0,\eta)}]-[u^{(0,0)}]).\] Consequently,  the set $\left\lbrace [u^{(0,\xi)}]-[u^{(0,0)}] | \xi \in \mathcal{A}(V) \right\rbrace$ is total in $\mathcal{H}(\mathcal{U}_\alpha)$. Hence, the prescribed map is a unitary, and $Ind(\alpha)=$dim~$\mathcal{A}(V)$. \hfill $\Box$

  \begin{rmrk}
  It is not clear to the authors how to meaningfully define the notion of index when $P$ is an arbitrary Lie semigroup. 
  \end{rmrk}

\section{Computation of additive cocycles}

In this section, we compute the units of the CCR flows associated to $P$-spaces, or equivalently, the additive cocycles of isometric representations associated to $P$-spaces. The case of a cone was treated in \cite{Anbu_Sundar}. 
We show in this section that, for unimodular groups, the isometric representation associated to a $P$-space $A$ has a non-zero additive cocycle if and only if the boundary of $A$ is compact.

Let $G$ be a connected, unimodular Lie group with Haar measure $\mu$, and let $P$ be a Lie semigroup of $G$ containing the identity $e$. Assume that $\Omega:=Int(P)$ is dense in $P$ and $PP^{-1}=P^{-1}P=G$.  Let $L(P)$ be the Lie-wedge of $P$.   By Corollary 3.11 of \cite{Neeb_Hilgert}, $\Omega$ is contained in the semigroup generated by $exp(L(P))$. Recall the preorder on $G$. For $x,y \in G$, we say  $x \leq y$ if $x^{-1}y \in P$ and $x < y$ if $x^{-1}y \in \Omega$. 

For $x \in P$ and $a \in \Omega$,  let 
\begin{align*}
[e,x]:&=\left\lbrace z \in G | e \leq z \leq x \right\rbrace = \left\lbrace z \in P | z \leq x \right\rbrace \\
(e,a):&=\left\lbrace z \in G | e<z<a \right\rbrace = \left\lbrace z \in \Omega | z<a \right\rbrace.
\end{align*}

\begin{rmrk}\label{monotone}
Let us collect a few preliminary observations that we need in this remark. 
\begin{enumerate}
\item[(1)] There exists a sequence $(a_k) \in \Omega$ such that $a_k<a_{k+1}$ and $(a_k)$ is cofinal in $G$. To see this, note that  $G=\Omega \Omega^{-1}= \displaystyle \bigcup_{b \in \Omega} \Omega b^{-1}$, a union of open sets. However, $G$ is second-countable, hence it can be expressed as a countable union of such open sets, $G= \bigcup_{n \geq 1} \Omega b_n^{-1}$. Then, $b_1^{-1}$ is contained in $\Omega b_{n_1}^{-1}$, for some $n_1$, and $b_1 < b_{n_1}$. Choose $x \in b_2 \Omega \cap b_{n_1} \Omega $; there exists $n_2$ such that $x^{-1} \in \Omega b_{n_2}^{-1}$. So we have, $b_2, b_{n_1} <x$ and $x<b_{n_2}$. Hence $b_2 < b_{n_2}$ and $b_{n_1} < b_{n_2}$. Repeating this argument for $ b_k$ and $b_{n_{k-1}}$,  we get a strictly increasing sequence $(a_k = b_{n_k})$ such that $\Omega b_k^{-1} \subset \Omega a_k^{-1}$ and $(a_k)$ is a cofinal sequence in $\Omega$. 

\item[(2)]
A continuous path $\gamma:[0,1] \to G$ is said to be monotone if whenever $s \leq t$, $\gamma(s) \leq \gamma(t)$.  We claim that any interior point in $P$ can be joined to the identity element $e$ by a continuous monotone path in $P$. To see this, let $a \in \Omega$ be given. Write  $a=exp X_{1}...exp X_{n}$ for some $X_{1},...,X_{n} \in L(P)$. Let $\gamma_{i}:[0,1]\longrightarrow P$ be defined by $\gamma_{i}(t)=exp (tX_{i})$, for $i=1,...,n$. Then, each $\gamma_{i}$ is a  monotone path in $P$ between $e$ and $expX_{i}$. Note that $\exp X_{1}\gamma_{2}$ is a monotone path joining $exp X_{1}$ and $exp X_{1}exp X_{2}$. Similarly, $expX_{1}expX_{2}\gamma_{3}$ is a monotone path between $expX_{1}expX_{2}$ and $expX_{1}expX_{2}expX_{3}$, and so on. On concatenating the paths $\gamma_{1}$, $expX_{1}\gamma_{2}$, ..., $expX_{1}...expX_{n-1}\gamma_{n}$, we  get a monotone path from $e$ to $a$ inside $P$.
\item[(3)] Note that both $e$ and $a$ belong to $\overline{(e,a)}$, for $a$ in $\Omega$. If $W$ is any open set containing $a$, then $W \cap \Omega$ is also an open set containing $a$, so $W \cap \Omega \cap a \Omega^{-1}$ will be non-empty since $a$ is a limit point of $a \Omega^{-1}$. Thus, $a \in \overline{(e,a)}$. Also, $W \cap \Omega \cap (\Omega^{-1} a)$ is non-empty by the same logic, so there exists a sequence $(t_n)$ in $\Omega \cap( \Omega^{-1} a)$ converging to $a$. Thus, $a t_n ^{-1}$ is a sequence in $\Omega \cap a \Omega^{-1}$ convering to $e$. This shows that $e$ belongs to $\overline{(e,a)}$.
\end{enumerate}
\end{rmrk}

Let $A$ be a $P$-space which will be fixed for the rest of this section. The interior of $A$ and the boundary of $A$ will be denoted $Int(A)$ and $\partial A$ respectively. Note that $\Omega A \subset Int(A)$ and $P Int(A) \subset Int(A)$. 
Consider the isometric representation $V=\left\lbrace V_{x} \right\rbrace _{x \in P}$ on $L^{2}(A,\mu)$, defined by \\
\begin{equation}
V_{x}f(y) = \begin{cases}
f(x^{-1}y) \quad \textrm{if~} x^{-1}y \in A,\\
0 \quad \textrm{otherwise}.
\end{cases}
\end{equation}
Here, as can be easily seen, $Ker(V_x^*) = L^2(A \backslash xA)$.

\begin{ppsn} \label{additive cocycle}

Suppose $\mu(A\backslash xA)<\infty$ for all $x \in P$. Then, for every complex number $\lambda$, $\{\lambda(1_{A\backslash xA})\}_{x \in P}$ is an additive cocycle of $V$. 

	Let $\xi: P \to L^{2}(A)$ be a map. If $\xi$ is an  additive cocycle of $V$, then there exists a complex number $\lambda$ such that, for $x \in P$, $\xi_{x}= \lambda (1_{A \backslash xA})$. 
		Consequently, $V$ has non-zero additive cocycles iff $\mu(A \backslash xA) < \infty$, for all $x \in P$.
\end{ppsn}
	\textit{Proof:} By considering a right translate of $A$, we can without loss of generality assume that $e \in A$, and consequently $P \subset A$. 	Fix a scalar $\lambda \in \mathbb{C}$. Let us assume that $A \backslash xA$ has finite measure for all $x \in P$. For $x \in P$, let $\xi_x:=\lambda 1_{A\backslash xA}$. 
Then $\xi_x \in Ker(V_x^*)$. Observe that for $x,y \in P$, $(A \backslash xyA) = (A \backslash xA) \coprod x(A \backslash yA)$. Hence, for $x,y \in P$, we have $\xi_{xy} = \xi_x + V_x \xi_y$. 

Let $(x_n)$ be a sequence in $P$ such that $(x_n)$ converges to $x$.  Since $A$ is closed, it follows that $1_{A\backslash x_nA}(y) \to 1_{A\backslash xA}(y)$ whenever $y \notin x \partial A$. Thanks to Lemma II.12 of \cite{Hilgert_Neeb}, it follows that $x\partial A$ has measure zero. Therefore, $(1_{A\backslash x_nA}) \to 1_{A \backslash xA}$ a.e.  
Consider an increasing cofinal sequence, say $(a_n)$, in $\Omega$. Then, $x < a_N$, for some $N$; $x^{-1}a_N \in \Omega$. Thus,  there exists $k$ such that $x_n^{-1}a_N \in \Omega$, for all $n \geq k$. For $n \geq k$, $(A \backslash x_n A) \subset (A \backslash a_NA)$ and hence, $1_{A \backslash x_n A}$ is bounded above by $1_{A \backslash a_NA}$.  By DCT, the sequence $(1_{A \backslash x_n A})$ converges to $1_{A \backslash xA}$ in $L^2(A ,\mu)$. Thus, the map $P \ni x \mapsto \xi_x \in L^2(A, \mu)$ is continuous. Therefore, $\xi= \left\lbrace \lambda(1_{A \backslash xA}) \right\rbrace _{x \in P}$ is an additive cocycle of $V$.

 	Now, let $\xi=\left\lbrace \xi_x \right\rbrace_{x \in P}$ be an additive cocycle of $V$. Each $\xi_x$ belongs to $L^2 (A \backslash xA)$. Let $(a_n)$ be an increasing cofinal sequence in $\Omega$ and set $a_0=e$. Then, $A=\displaystyle \coprod_{n \geq 1}( a_{n-1}A \backslash a_n A)$. Define  $f: A \longrightarrow \mathbb{C}$  by, 
\[ f(x):= \xi_{a_n}(x), \quad x \in a_{n-1} A \backslash a_n A. \]
The map $f$ is well-defined, and is also measurable by definition.
\begin{enumerate}
	\item[(1)] For all $n$, $f = \xi_{a_n}$, a.e. on $A \backslash a_n A$.
	Fix $n$. Note that $A \backslash a_n A = \bigcup_{1 \leq k \leq n} ( a_{k-1} A \backslash a_k A)$. For $m<n$ and $x \in a_{m-1}A \backslash a_m A$, $a_m^{-1}x \notin A$, so $V_{a_m} \xi_{a_m^{-1}a_n}(x)=0$. Thus, for $m<n$
	\[ \xi_{a_n}(x) = \xi_{a_m}(x) + V_{a_m} \xi_{a_m^{-1}a_n}(x) = \xi_{a_m}(x) = f(x)\]
	for almost all $x \in a_{m-1}A \backslash a_m A$. Hence, $f(x) = \xi_{a_n}(x)$ a.e. on $A \backslash a_n A$, for all $n$.
	\item[(2)]  Let $a \in P$ be given. Then  $f(x) = \xi_a(x)$, a.e. on $A \backslash aA$.
	Fix $a$ in $P$. Choose $n$ for which $a < a_n$, so $(A \backslash aA) \subset (A \backslash a_n A)$. But, for $x$ in $A \backslash aA$, $a^{-1}x \notin A$, so $V_a \xi_{a^{-1}a_n}(x)=0$. Hence, for almost every $x \in (A \backslash aA)$, 
	\[ \xi_{a_n}(x)= \xi_a(x) + V_a \xi_{a^{-1}a_n}(x)= \xi_a(x).\]
	Thus, $f (x)= \xi_a(x)$ for almost all $x \in A \backslash aA$.
	\item[(3)] For every $a \in P$, $f(ax)=f(x)$, a.e. on $A$.
	Fix $a$ in $P$. Note that, for all $x \in A$, $\xi_a (ax)=0$. Thus, for each $n$, 
	\[\xi_{aa_n}(ax) = \xi_a(ax) + V_a \xi_{a_n}(ax) = \xi_{a_n}(x)\]
	for a.e. $x \in A$. But by $(2)$, $\xi_{a_n}(x) = f(x)$ and $\xi_{aa_n}(ax)=f(ax)$ for a.e $x$ in $A \backslash a_n A$. Since $A= \bigcup_n (A \backslash a_n A)$, $f(ax)=f(x)$ for almost every $x \in A$.
\end{enumerate}
	
We define another order on $G$ as follows. For $g, h\in G$, we say $g \leq_r h$  if $g^{-1}h \in P^{-1}$, or equivalently $h^{-1}g \in P$. With this order, $\Omega^{-1}$ has an increasing cofinal sequence, say $(b_n)$, and we set $b_0=e$. For $g \in G$, $g \leq_r b_n$, for some $n$.  Hence, $g \in b_n P \subset b_n A$. Therefore, $G=\bigcup_{n \geq 0}b_nA$.  Also, for $m <n$, we have $b_m A \subset b_n A$.  We extend $f$ to $G$ by:
\begin{align*}
	\widetilde{f}|_{(b_0 A =A)}(x)&=f(x),\\
	\widetilde{f}|_{(b_n A \backslash b_{n-1}A)}(x)&= f(b_n^{-1}x).
\end{align*}
Then, $\widetilde{f}$ is well-defined, measurable, and is an extension of $f$ to $G$.
\begin{enumerate}
	\item[(4)] For all $n$, $\widetilde{f}(x) = f(b_n ^{-1} x)$, for almost all $x \in b_n A$. 	Fix a natural number $n$. Observe that $b_n A = \displaystyle A \coprod_{1 \leq k \leq n} (b_k A \backslash b_{k-1}A)$. By (3),
	\[ f(b_n^{-1}x) = f(b_n^{-1} b_k b_k^{-1}x) = f(b_k^{-1}x) = \widetilde{f}(x)\]
	for almost every $x \in (b_k A \backslash b_{k-1} A)$ and for all $1 \leq k \leq n$. Thus, $\widetilde{f}|_{b_n A}(x)=f(b_n^{-1}x)$ for almost every $x \in b_n A$. 
	\item[(5)] Let $b \in P^{-1}$ be given. We claim that  $\widetilde{f}(x)=f(b^{-1}x)$ for almost every $x \in bA$. 
	Since $(b_n)$ is cofinal, it follows that for some $n$, $b \leq_r b_n$, i.e. $b_n^{-1}b \in P$. Then, by (3) and (4), for almost every $x \in bA \subset b_nA$, \[\widetilde{f}(x)=f(b_n^{-1}x)=f(b_n^{-1}bb^{-1}x)=f(b^{-1}x).\]
\end{enumerate}
     Fix $b \in P^{-1}$. For almost all $x \in b_nA$,
     \[\widetilde{f}|_{bb_nA}(bx) = f((bb_n)^{-1}bx) = f(b_n^{-1}x) = \widetilde{f}(x), \]
     for all $n$. But $G = \displaystyle \bigcup_n b_nA$, and hence $\widetilde{f}(bx) = \widetilde{f}(x)$ for almost all $x \in G$. Since $ P^{-1}$ generates $G$, it follows that for every $g \in G$, 
     $
     \widetilde{f}(gx)=\widetilde{f}(x)$
     for almost all $x \in G$. But the left translation action of $G$ on $G$ is ergodic. Therefore, there exists a complex number $\lambda$ such that $\widetilde{f}=\lambda$ a.e.  By (2),  for $x$ in $P$, $\xi_x = f|_{(A \backslash xA)} = \lambda 1_{(A \backslash xA)}$. The proof is now complete.  \hfill $\Box$

 \begin{lmma} \label{onto}
 Let $a \in \Omega$ be given.  We have the following. 
	 \begin{itemize}
		\item[(1)]  The set $(e,a)$ is path-connected.
		\item[(2)]   The map $ \partial A \times (e,a) \ni (x,b) \to bx \in  Int A \smallsetminus aA$ is a surjection. 
	\end{itemize}
\end{lmma}
 \textit{Proof.} Let $b,c \in (e,a)$ be given. Choose  a sequence $(s_{n}) \in \Omega$ such that $ (s_{n})$ converges to $ e$. Since $b \in \Omega$ and $(s_{n}^{-1}b)\to b \in \Omega$, we have  $s_n^{-1}b \in \Omega$ for large $n$. Similarly $s_n^{-1}c \in \Omega$ for large $n$. Choose $N$ such that $s_N^{-1}b \in \Omega$ and $s_{N}^{-1}c \in \Omega$. Then,   $e<s_{N}<b<a$.  By Remark \ref{monotone}, there exists a monotone path  $\sigma:[0,1]\longrightarrow P$  joining $e$ and $s_{N}^{-1}b$. Then $s_{N}\sigma$ is a monotone path in $(e,a)$ joining $s_{N}$ and $b$. Similarly,  there exists a path in $(e,a)$ joining $s_{N}$ and $c$. Thus, we have a path in $(e,a)$ joining $b$ and $c$. Therefore, $(e,a)$ is path-connected. This proves (1). 

 Next, to see that the map defined in $(2)$ is well defined,  let $b \in (e,a) \subset \Omega$ and $x \in \partial(A) \subset A$ be given.  Then, $bx \in Int(A)$. If $bx \in aA$, then $bx=ak$ for some $k \in A$ and $x=b^{-1}ak$. Since $b \in (e,a)$, we have  $b^{-1}a \in \Omega$. But $\Omega A \subset Int(A)$. Therefore, $x \in Int(A)$, which is a contradiction. Hence $bx \in Int (A) \smallsetminus aA$ and our map is well-defined which we call as $\psi$. 
 
To prove the surjectivity of $\psi$, let  $b \in Int (A) \smallsetminus aA$ be given.   
Consider the map $\phi : (e,a) \longrightarrow G$ defined by $\phi(s)=s^{-1}b$.  We claim that $\phi((e,a)) \cap A^{c}$ is non-empty. Thanks to Remark \ref{monotone}, there exists a sequence $(s_n)$ in $(e,a)$ such that $s_n \to a$. Note that $s_n^{-1}b \to a^{-1}b$, $a^{-1}b \in A^{c}$ and $A^{c}$ is open. Thus, there exists $N \geq 1$ such that $s_{n}^{-1}b \in A^{c}$ for all $n\geq N$. This proves that $\phi((e,a)) \cap A^{c}$ is non-empty.

Let $(t_{n})$ be a sequence in $(e,a)$ that converges to $e$. Since $b \in Int(A)$, $t_{k}^{-1}b \in Int (A)$, for large $k$. Hence, $\phi ((e,a))\cap Int (A)$ is non-empty.  However, $\phi((e,a))$ is connected as $(e,a)$ is connected. Therefore,  $\phi ((e,a)) \cap \partial A$ is non-empty. This implies that there exists $s \in (e,a)$ such that $s^{-1}b \in \partial A$ and $\psi(s^{-1}b,s)=ss^{-1}b=b$. Consequently, the map $\psi$ is onto. This completes the proof.  \hfill $\Box$

\begin{ppsn}
\label{finite vs compact}
The following are equivalent. 
\begin{enumerate}
\item[(1)] For every $x \in P$, $A \backslash xA$ has finite measure. 
\item[(2)] For every $a \in \Omega$, $A\backslash aA$ has finite measure. 
\item[(3)] The boundary $\partial A$ is compact. 
\end{enumerate}
\end{ppsn}
\textit{Proof.}   Assume that $(2)$ holds. Suppose $\partial A$ is not compact. Fix $a \in \Omega$. Choose a compact set $E \subset (e,a)$ with non-empty interior. Since $\partial A$ is not compact, there exists a sequence $(x_{n}) \in \partial A$ which has no convergent subsequence. We claim that there exists $k_1 >1$ such that $Ex_{k_1} \cap Ex_1=\emptyset$. 
Because if not, then $Ex_{n} \cap Ex_{1} \neq \emptyset$, for all $n$. Choose $y_{n} \in Ex_{n} \cap Ex_{1}$. Then, $y_{n}=h_{n}x_{n}=g_{n}x_{1}$, for $h_{n}, g_{n} \in E$, $x_n x_1^{-1}=g_{n}h_{n}^{-1}$. But $g_{n}h_{n}^{-1} \in EE^{-1}$, which is compact. Thus, $(x_{n}x_{1}^{-1})$, and consequently $(x_{n})$, will have a convergent subsequence, which is a contradiction. Hence, there exists $k_{1}$ such that $Ex_{k_{1}} \cap Ex_{1}$ is the null set. Let \[I_{k_{1}}= \left\lbrace k \in \mathbb{N} | k \geq k_{1}, Ex_{k} \cap Ex_{k_{1}} = \emptyset, Ex_{k} \cap Ex_{1} = \emptyset \right\rbrace .\] By a very similar arguement, we see that $I_{k_{1}}$ is non-empty. So, there exists $ k_{2} >k_1 >1$ such that $ Ex_{1} \cap Ex_{k_{2}} = \emptyset$ and $Ex_{k_{1}} \cap Ex_{k_{2}} = \emptyset$. We continue this process to get a subsequence $(x_{n_{k}})$ such that $Ex_{n_{k}} \cap Ex_{n_{m}} = \emptyset$, for all $k \neq m$.  Since $G$ is unimodular, $\mu(Ex_{n_k})=\mu(E)$.  As verified in Lemma \ref{onto}, $Ex_{n_{k}}  \subset Int~A \smallsetminus aA$ for all $k$. Consequently, $Int(A)\backslash aA$ contains the disjoint union $\displaystyle \coprod_{k=1}Ex_{n_k}$ and the latter set has infinite measure. Therefore,  $\mu(Int A \smallsetminus aA) = \infty$. This completes the proof of the implication $(2) \implies (3)$.

Suppose that $(3)$ holds.  Set $H:=P \cap P^{-1}$. We claim that $H$ is compact. Since $PA \subset A$, it follows that  for $h \in H$, $hA \subset A$ and $h^{-1}A \subset A$. In other words, $hA=A$ for every $h \in H$. Consequently,  $h \partial A=\partial A$ for every $h \in H$. Fix $x_0 \in \partial A$. Then the map 
\[
 H \ni h \to hx_0 \in \partial A\]
 is a topological embedding. Since $\partial A$ is compact, it follows that $H$ is compact. This proves the claim. 

Let $\widetilde{G}=G / H$ be the homogeneous space of left cosets of $H$. For $x \in G$, we denote the left coset $xH$ by $\widetilde{x}$.  The map $G \ni x \to \widetilde{x} \in \widetilde{G}$ will be denoted by $\pi$.  The preorder $\leq$ descends to a closed partial order on $\widetilde{G}$. That is, for $x,y\in G$, $\widetilde{x} \leq \widetilde{y}$ if $ x \leq y$.  It is easily verifiable using Remark \ref{monotone} that $\widetilde{G}$ has the chain approximation property defined in  Page 116, Chapter 4 of \cite{Neeb_Hilgert}. Making use of  Prop. 4.4 of \cite{Neeb_Hilgert},  choose an open set  $\widetilde{U}$ containing $\widetilde{e}$ such that for all $ \widetilde{a} \in \widetilde{U}$, $[\widetilde{e}, \widetilde{a}] $ is compact whenever $\widetilde{a} \geq \widetilde{e}$. . Let $U:=\pi^{-1}(\widetilde{U})$. Note that $U$ contains $e$. 
Since $H$ is compact, it follows that $\pi$ is proper, i.e. the inverse image of a compact set is compact. Therefore, for $a \in U$, $[e,a]=\pi^{-1}([\widetilde{e},\widetilde{a}])$ is compact. 

\textit{Claim:} For $a \in U \cap P$,  $Int(A)\smallsetminus aA$ has finite measure.

Let $a \in U \cap P$ be given. 

\textbf{Case 1:} Suppose $a \in \Omega$. Consider the map $\psi:\partial(A) \times [e,a]  \to vx \in G$ defined by $\psi(x,v)=vx$. Then $\psi$ is  a continuous map and has  compact image. Thanks to Lemma \ref{onto}, the image of $\psi$ contains $Int(A)\backslash aA$. Hence, $ Int(A) \smallsetminus aA$ has finite measure.

\textbf{Case 2:} Suppose $a \in \partial P$. Since $e \in U$, there exists a sequence $(s_n)$ in $U \cap \Omega$ that converges to $e$. Thus, $(as_n)$ converges to $a \in U$. Hence, there exists $N$ for which $as_N \in U$. Clearly $as_N \in \Omega$. By Case 1, $Int(A)\backslash as_N A$ has finite measure. Note that 
\[
Int(A)\backslash as_{N} A=(Int(A) \backslash aA) \coprod a(A \backslash s_N A).\]
Hence $Int(A)\backslash aA$ has finite measure. 

This proves our claim. 

Note that by Lemma II.12 of \cite{Hilgert_Neeb}, $\partial A$ has measure zero. Therefore, for $a \in U \cap P$, $A\backslash aA$ has finite measure. This is because, up to a set of measure zero,  $A \backslash aA=Int(A)\backslash aA$. 

Let $a \in \Omega$ be given. Write $ a=expX_{1}expX_{2}...expX_{n}$, for some $X_{1},...X_{n} \in L(P)$. Since $U$ is  open and contains  $e$, there exists a natural number $N$ such that for every  $i$, $exp(\frac{X_i}{N}) \in U$. By what we have proved above, 
$\mu(A\backslash exp(\frac{X_i}{N}))<\infty$. 

But $A \backslash exp(X_i)A$ is the disjoint union of the  sets $A \backslash exp(\frac{X_i}{N})A$, $exp(\frac{X_i}{N})(A \backslash exp(\frac{X_i}{N})A), \cdots $ $exp(\frac{(N-1)X_i}{N})(A \backslash exp(\frac{X_i}{N})A)$ each having finite measure. Therefore, for each $i$,  $A \backslash exp(X_i)A$ has finite measure.  Note that  $A \backslash aA$ is the disjoint union of the  sets $A \backslash exp(X_1)A$, $exp(X_1)(A \backslash exp(X_2)A)$, $\cdots$, $exp(X_1)exp(X_2)\cdots exp(X_{n-1})(A \backslash exp(X_n)A)$.  Hence $A \backslash aA$ has finite measure for every $a \in \Omega$. 

Let $x \in P$ be given. Choose $s \in \Omega$. Then $A \backslash xA \subset A \backslash xs A$. But $xs \in \Omega$ and $A \backslash xs A$ has finite measure. Therefore, $A \backslash xA$ has finite measure for every $x \in P$. This completes the proof of the implication $(3) \implies (2)$. \hfill $\Box$
 
\begin{crlre} \label{sec 3 result}
	The isometric representation $V$ has a non-zero additive cocycle iff $\partial A$ is compact. 
\end{crlre}
	\textit{Proof:} This is immediate from Prop. \ref{additive cocycle} and Prop. \ref{finite vs compact}. \hfill $\Box$

\section{Is the boundary compact ?}

In this section, we discuss whether the boundary of a $P$-space is compact or not. We derive a necessary condition for the boundary to be compact which we prove is also sufficient in the abelian case. 
Recall the Malcev-Iwasawa theorem (\cite{Borel}): Let $G$ be a connected Lie group and let $K$ be a maximal compact subgroup of $G$. Then $G$ is topologically homeomorphic to $K \times \mathbb{R}^n$ for some $n$, and $G/K$ is topologically homeomorphic to $\mathbb{R}^n$. 

\begin{thm} \label{G/K}
	Let $G$ be a connected Lie group and let $K$ be a maximal compact subgroup. Suppose $P$ is a Lie semigroup such that $PP^{-1}=P^{-1}P=G$, and $A$ is a $P$-space. If the boundary of $A$ is compact, then $dim ~ (G/K) = 1$. 
\end{thm}
\textit{Proof:} There is no loss of generality in assuming  that $A$ contains $P$. By the Malcev-Iwasawa Theorem, we have $G= K \times \mathbb{R}^n$, up to a homeomorphism, for some $n$. Let the boundary of $A$, $\partial A$, be compact. Then $\partial A \subset K \times B[0,R]$, for some large $R>0$, where $B[0,R]$ is the closed ball of radius $R$ in $\mathbb{R}^n$. Now,
\[ G \backslash (K \times B[0,R]) = ( Int(A) \backslash (K \times B[0,R])) \cup (A^c \backslash (K \times B[0,R])).\]

Observe that $Int(A)$ contains $\Omega$, which does not have compact closure. Note that $P^{-1}A^c \subset A^c$. Therefore,  $A^{c}$ contains a translate of $\Omega^{-1}$ which again does not have compact closure. So, both $Int(A) \backslash (K \times B[0,R])$ and $A^c \backslash (K \times B[0,R])$ are non-empty open sets, which do not intersect. Consequently, $G \backslash (K \times B[0,R])$ is not a  connected set and we have $n=dim ~ (G / K) = 1$. Hence the result. \hfill $\Box$

We do not know whether the necessary condition of Theorem \ref{G/K} is also sufficient.  However, we show that it is sufficient when $G$ is abelian.
For the rest of this section, let $G$ be a non-compact abelian Lie group. Then $G$ is of the form $G=\mathbb{R}^{d}\times\mathbb{T}^{r}$, for some non-negative integers $d \geq 1$ and $r \geq 0$. Let $P$ be a Lie semigroup of $G$ with dense interior $\Omega$.  Let $A$ be a $P$-space. We may assume that $A$ contains  $P$.

 Let $\exp:\mathbb{R}^{d}\times\mathbb{R}^{r}\longrightarrow\mathbb{R}^{d}\times\mathbb{T}^r$ be the exponential map, i.e.
\[ \exp ( x_{1},x_{2},...,x_{d},y_{1},y_{2},...y_{r} ) = ( x_{1},x_{2},...,x_{d},e^{2\pi i y_{1}},e^{2\pi i y_{2}},...,e^{2\pi i y_{r}}) .\]
The map $\exp$ is a homomorphism and its kernel is $(0,0,...,0)\times\mathbb{Z}^{r}$.
The Lie wedge of $P$ is 
\[L(P)=\left\lbrace( x_{1},...,x_{d},y_{1},...y_{r})\in\mathbb{R}^{d+r}: \exp(t( x_{1},x_{2},...,x_{d},y_{1},y_{2},...y_{r}))\in P,\forall t\geq0\right\rbrace.\]
Then $L(P)$ is a closed convex cone in $\mathbb{R}^{d+r}$ which contains the origin, has a dense interior, and is spanning in $\mathbb{R}^{d+r}$. Since $G$ is abelian, $P$ is a \textit{divisible} subsemigroup of $G$.  (Recall that we say $P$ is divisible, if for every $a \in P$ and $N \geq 1$, there exists $x \in P$ such that $x^{N}=a$). By Theorem V.6.5 of \cite{Hilgert_Hofmann_Lawson}, $P=\exp(L(P))$. 

We define $\pi:\mathbb{R}^{d}\times\mathbb{R}^{r}\longrightarrow\mathbb{R}^{d}$ and $\widetilde{\pi}: \mathbb{R}^d \times \mathbb{T}^r \longrightarrow \mathbb{R}^d$ as the projection maps,
\begin{align*}
\pi(x_{1},x_{2},...,x_{d},y_{1},...,y_{r}) &:=(x_{1},x_{2},...,x_{d}), \\
\widetilde{\pi}(x_1,x_2,...,x_d,w_1,w_2,\cdots,w_r) &:=(x_1,x_2,...,x_d).
\end{align*}
Note that $\widetilde{\pi} \circ \exp =\pi $. 
The map $\widetilde{\pi}$ is closed, since $\mathbb{T}^r$ is compact. So, $\widetilde{\pi}(P)=\pi(L(P))$ is closed in $\mathbb{R}^d$. Hence, the set $\pi(L(P))$ is a closed convex cone that spans $\mathbb{R}^{d}$. Set $P_1 = \pi(L(P)) = \widetilde{\pi}(P)$ and $A_1 = \widetilde{\pi}(A)$. Then, $A_{1}$ is a $P_{1}$-space and contains $P_{1}$.

For $x=(x_{1},...,x_{d}) \in P_1$ and for $a=(a_{1},...,a_{d}) \in A_{1}$, define 
\begin{align*} S(x)&:=\left\lbrace(y_{1},...,y_{r})\in\mathbb{R}^{r}|(x_{1},...,x_{d},y_{1},...,y_{r})\in L(P)\right\rbrace, \\
T(a)&:=\left\lbrace(b_{1},...,b_{r})\in\mathbb{R}^{r}| \exp((a_{1},...,a_{d},b_{1},...,b_{r}))\in A\right\rbrace.
\end{align*}
Note that $L(P)=\displaystyle \bigcup_{x\in P_{1}} (\left\lbrace x\right\rbrace \times S(x))$, and $A=\displaystyle \bigcup_{a\in A_{1}} \exp(\left\lbrace a \right\rbrace \times T(a))$.

\begin{lmma} \label{observations}
With the above notation, the following hold.
\begin{itemize}
	\item[(1)] For  $t\geq0$ and $x\in P_{1}$, $tS(x)\subset S(tx)$.
	\item[(2)] For  $a\in A_{1}$ and $x\in P_{1}$, $T(a)+S(x)\subset T(a+x)$. 
	\item[(3)]  The interior $Int (L(P))$ is contained in $\displaystyle \bigcup_{x\in \pi(Int  L(P))} (\left\lbrace x\right\rbrace \times Int ~S(x))$.
\end{itemize}
\end{lmma}
\textit{Proof:} Let $x=(x_{1},...,x_{d})$ in $P_{1}$ and $(y_{1},...,y_{r})$ in $S(x)$. Since $L(P)$ is a cone,  it follows that for $t \geq 0$, $t(x_{1},...,x_{d},y_{1},...,y_{r}) \in L(P)$. Thus, $t(y_{1},...,y_{r})$ belongs to $S(tx)$, for all $t\geq 0$, and we get $tS(x)\subset S(tx)$. This proves (1).

To prove (2), let $a=(a_{1},...,a_{d})\in A_{1}$ and $x=(x_{1},...,x_{d})\in P_{1}$ be given. Let $(b_{1},...,b_{r})\in T(a)$ and $(y_{1},...,y_{r})\in S(x)$ be given. Then, $(a_{1},...a_{d},e^{2\pi i b_{1}},...,e^{2\pi i b_{r}})\in A$ and $(x_{1},...,x_{d},e^{2\pi  iy_{1}},...,e^{2\pi i y_{r}})\in P$. Since $A+P \subset A$, it follows that   \[(a_{1}+x_{1},...,a_{d}+x_{d},e^{2\pi i (b_{1}+y_{1})},...,e^{2\pi i (b_{r}+y_{r})})\in A.\] Thus, $(b_{1}+y_{1},...,b_{r}+y_{r})=(b_{1},...,b_{r})+(y_{1},...,y_{r}) \in T(a+x)$.
  Consequently, the inclusion $T(a)+S(x)\subset T(a+x)$ holds. 
  
For proving (3), let $x=(x_{1},...,x_{d})$ and $y=(y_{1},...,y_{r})$ be such that $(x;y)\in Int ~L(P)$. Then there exists $R>0$ such that $(x;y) \in B(x,R) \times B(y,R) \subset Int ~ L(P)$ , in $\mathbb{R}^{d+r}$, where $B(x,R)$ and $B(y,R)$ are open balls of radius $R$ centered at $x$ and $y$ in $\mathbb{R}^d$ and $\mathbb{R}^r$ respectively. So, $x \in \pi(B(x,R) \times B(y,R)) \subset \pi(Int ~ L(P))$. Also, for $z \in B(y,R)$, $(x;z) \in B(x,R) \times B(y,R) \subset L(P)$. Thus, $y \in B(y,R) \subset S(x)$, i.e. $y \in Int ~ S(x)$. This proves (3). \hfill $\Box$

\begin{lmma} \label{tube}
	Let $x \in \pi(Int(L(P)))$ be given. Then there  there exists $t_x > 0$ such that for all $t \geq t_x$ and $a \in A_1$, $\exp(\left\lbrace a+tx\right\rbrace \times T(a+tx))=\left\lbrace a+tx\right\rbrace \times \mathbb{T}^{r}$.
\end{lmma}
\textit{Proof:} Let $x \in \pi(Int ~ L(P))$. From (3) of Lemma \ref{observations}, we see that $Int ~S(x)$ is non-empty. Hence, $S(x)$ contains a r-dimensional hypercube of side $l$, say. Using (1) from Lemma \ref{observations}, we can thus obtain a $t_{x}>0$ such that  $S(tx)$ contains a r-dimensional hypercube of side $1$ in $\mathbb{R}^r$, for all $t\geq t_{x}$.

Now, let $a \in A_{1}$ and $t\geq t_{x}$. We can write $T(a)+S(tx)=\bigcup_{b\in T(a)}(b+S(tx))$. Since $S(tx)$ contains a r-dimensional hypercube of side $1$,  each set $b+tS(x)$ contains a r-dimensional hypercube of side $1$ in $\mathbb{R}^r$. But $T(a)+S(tx)$ is the union of such sets, hence it also contains a r-dimensional hypercube of side $1$. Using (2) from Lemma \ref{observations}, we get that $T(a+tx)$ conatins a r-dimensional hypercube of side $1$ in $\mathbb{R}^r$. Therefore, under the exponential map,
\[
\exp(\left\lbrace a+tx\right\rbrace \times T(a+tx))=\left\lbrace a+tx\right\rbrace \times \mathbb{T}^{r}, \quad \forall ~ t\geq t_{x}. \]
This completes the proof. \hfill $\Box$

Fix $x_{0} \in \pi(Int ~L(P))$. Let  $t_0:=t_{x_{0}}$  be as in Lemma \ref{observations}. Thus, for all $a \in A_{1}$,
\[ \exp(\left\lbrace a+t_{0}x_{0}\right\rbrace \times T(a+t_{0}x_{0}))=\left\lbrace a+t_{0}x_{0}\right\rbrace \times \mathbb{T}^{r}.\]

\begin{rmrk}
	Suppose $A_{1}=\mathbb{R}^{d}$. Then $a-t_{0}x_{0}$ $\in$ $A_{1}$, for all $a \in \mathbb{R}^{d}$. Thus, 
	\[\exp(\left\lbrace a\right\rbrace \times T(a))=\exp(\left\lbrace a-t_{0}x_{0}+t_{0}x_{0}\right\rbrace \times T(a-t_{0}x_{0}+t_{0}x_{0}))=\left\lbrace a\right\rbrace \times \mathbb{T}^{r}.\]
	Since $A=\bigcup_{b\in A_{1}}\exp(\left\lbrace b \right\rbrace \times T(b))$, we  have $A=\mathbb{R}^{d} \times \mathbb{T}^{r}$ which is a contradiction. 
	
	Similarly, if $P_1=\widetilde{\pi}(P)=\mathbb{R}^{d}$, then  $A_{1}$ being a $P_{1}$-space and containing $P_{1}$, we  have $A_{1}=\mathbb{R}^{d}$ and thus $A=\mathbb{R}^{d} \times \mathbb{T}^{r}$ which is again a contradiction.
	Therefore, both $P_1$ and $A_1$ are proper subsets of $\mathbb{R}^{d}$.
\end{rmrk}

\begin{thm}
\label{when is boundary}
	The boundary of $A$, $\partial A$, is compact if and only if $d=1$.
\end{thm}

\textit{Proof:} Let $\partial A$ be compact. Since $G = \mathbb{R}^d \times \mathbb{T}^r$, the maximal compact subgroup of $G$ is $\mathbb{T}^r$. Hence, by Theorem \ref{G/K}, $d=1$.

To prove the converse, let us assume $d=1$. Then, $P_{1}=[ 0,\infty)$ or $P_{1}=(-\infty ,0]$. Consider the case when $P_{1}=[0, \infty)$. Choose $x_{0}$ from $\pi(Int ~L(P))$. Let $t_0:=t_{x_0}$ be as in Lemma \ref{tube}.  Set $x:=t_0x_0$. 

Since $A_{1}$ is not entire $\mathbb{R}$ and $A_1$ is a $P_1$-space, $A_{1}=[k, \infty)$ for some $k \leq 0$. By Lemma \ref{tube}, for all $a\in A_{1}$, $\exp(\left\lbrace a \right\rbrace \times T(a))= \left\lbrace a \right\rbrace \times \mathbb{T}^r$, whenever  $a\geq (k+x)$, and we get 
\[ [k+x, \infty) \times \mathbb{T}^{r} \subset A \subset [ k, \infty) \times \mathbb{T}^{r} .\]
Thus, $\partial A$ is contained in the compact set $[ k, k+x] \times \mathbb{T}^{r}$, and hence is compact.
The case $P_{1}=(-\infty, 0]$ is similar. \hfill $\Box$

\section{The CCR flows $\alpha^{(A,k)}$}
In this section, we provide another direct proof of the fact for a pure isometric representation $V$, the corresponding CCR flow $\alpha^{V}$ remembers the representation $V$. We prove this in the setting
of Lie semigroups. We then apply this to classify the CCR flows $\alpha^{(A,k)}$, or equivalently the corresponding isometric representations $V^{(A,k)}$.

Let $G$ be a connected Lie group. Let $P$ be a Lie semigroup of $G$ having a dense interior $\Omega$ such that $PP^{-1}=P^{-1}P=G$. 
Let us recall the product system of a CCR flow. Let $V: P \to B(\clh)$ be a pure isometric representation. Recall that $V$ is said to be \emph{pure} if $\displaystyle \bigcap_{x \in P}Ran(V_x)=\{0\}$. Let $\alpha^{V}:=\{\alpha_x\}_{x \in P}$ be the CCR flow associated to $V$. 
Denote the product system of $\alpha^{V}$ by $E:=\{E(x)\}_{x \in P}$. 

Then, for $x \in P$, we can identify $E(x)$ with $\Gamma(Ker(V_x^{*}))$, where $\Gamma(Ker(V_x^{*}))$ is the symmetric Fock space of $Ker(V_x^*)$. Moreover, the product on the exponential vectors is given by
\[
e(\xi)e(\eta)=e(\xi+V_x\eta)
\]
for $\xi \in Ker(V_x^{*})$ and $\eta \in Ker(V_y^{*})$.

\begin{thm}
\label{injectivity}
	Let $V:P\longrightarrow B(\mathcal{H})$ and $W:P \longrightarrow B(\mathcal{K})$ be pure isometric representations . The CCR flows associated to $V$ and $W$ are cocycle conjugate iff $V$ and $W$ are unitarily equivalent.
\end{thm}
\textit{Proof:} Let $\alpha$ and $\beta$ be the CCR flows associated to $V$ and $W$ respectively.
 Assume that $\alpha$ and $\beta$ are cocycle conjugate. Denote the product systems of $\alpha$ and $\beta$ by $\mathcal{E}^\alpha = \displaystyle \coprod_{x\in P} E^{\alpha}(x)$ and $\mathcal{E}^\beta = \displaystyle \coprod_{x \in P} E^\beta(x)$ respectively. Then, there exists a Borel isomorphism $\theta: \mathcal{E}^\alpha \longrightarrow \mathcal{E}^\beta$ such that $\theta(uv)= \theta(u) \theta(v)$ for $u \in E^\alpha(x)$ and $v \in E^\alpha(y)$ and
\[\theta_x : E^\alpha(x) \longrightarrow E^\beta(x), \quad \theta_x=\theta|_{E^\alpha(x)}, \] is a unitary for all $x\in P$. 

We claim that for $a \in \Omega$ and $\xi \in Ker(V_a^{*})$, there exist a unique non-zero complex number $c_{a,\xi}$ and a unique vector $\widetilde{\xi}_{a} \in Ker(W_a^{*})$ such that 
\[
\theta_{a}(e(\xi))=c_{a,\xi}e(\widetilde{\xi}_{a}).\]
The proof of the claim is well known if $P=[0,\infty)$. It follows, for instance, from Arveson's computation of the gauge  group of a $1$-parameter CCR flow. Let $X \in L(P)$ be given. Restricting the product systems to the ray $\{exp(tX): t \geq 0\}$ and applying this one parameter assertion, we see that given $\xi \in Ker(V_{exp(X)}^{*})$, there exist a non-zero complex number $c_{X,\xi}$ and a vector $\widetilde{\xi}_{X} \in Ker(W_{exp(X)}^{*})$ such that 
\[
\theta_{exp(X)}e(\xi)=c_{X,\xi}e(\widetilde{\xi}_{X}).\]

Fix $a \in \Omega$. Since $\Omega$ is contained in the semigroup generated by $\exp(L(P))$, there exist $X_1, X_2, ..., X_n$ in $L(P)$ such that $a=\exp(X_1)\exp(X_2)...\exp(X_n)$. We can expand $Ker(V_a^*)$ as 
\[ Ker( V_a^*) = Ker(V_{\exp(X_1)}^*) \oplus V_{\exp(X_1)}Ker(V_{\exp(X_2)}^*) \oplus \cdots \oplus V_{\exp(X_1) \exp(X_2)...\exp(X_{n-1})} Ker(V_{\exp(X_n)}^*).\]
Thus, for any $\xi \in Ker(V_a^*)$, there exist $\xi_i \in Ker( V_{\exp(X_i)}^*)$, $1 \leq i \leq n$ such that
\[\xi = \xi_1 + V_{\exp(X_1)} \xi_2 + ... +V_{\exp(X_1)...\exp(X_{n-1})} \xi_n.\]
Note that $e(\xi)=e(\xi_1)e(\xi_2)\cdots e(\xi_n)$. Since $\theta$ is multiplicative, it follows that 
\begin{align*}
\theta_a e(\xi) &= \theta_{exp(X_1)} e(\xi_1) \theta_{exp(X_2)} e(\xi_2)
... \theta_{exp(X_n)} e(\xi_n) \\
&= c_{X_1, \xi_1}e(\widetilde{\xi}_1) c_{X_2, \xi_2}e(\widetilde{\xi}_2) ...c_{X_n, \xi_n}e(\widetilde{\xi}_n) \\
&= c_{X_1, \xi_1} c_{X_2, \xi_2}  ...  c_{X_n, \xi_n}e( \widetilde{\xi}_1 + W_{\exp(X_1)} \widetilde{\xi}_2 + ... + W_{\exp(X_1)...\exp(X_{n-1})} \widetilde{\xi}_n) \\
&= c_{a,\xi} e (\widetilde{\xi}_a)
\end{align*}
where \[\widetilde{\xi}_a:=\widetilde{\xi}_{1}+W_{exp(X_1)}\widetilde{\xi}_{2}+\cdots +W_{exp(X_1)\cdots exp(X_n)}\widetilde{\xi}_n \in Ker(W_{a}^*)\] and $c_{a, \xi} := c_{X_1, \xi_1}  c_{X_2, \xi_2}  ...  c_{X_n, \xi_n}$. Uniqueness of $c_{a,\xi}$ and $\widetilde{\xi}_{a}$ is clear. This proves the claim. 

Now, for any $a, b \in \Omega$,
\begin{align*}
\theta_{ab} e(0) &= \theta_a e(0) \theta_b e(0) \\
&= c_{a, 0} e(\widetilde{0}_a) c_{b, 0} e (\widetilde{0}_b)\\
&= c_{a,0}c_{b,0} e(\widetilde{0}_{a} + W_a \widetilde{0}_b).
\end{align*}
Thus, for all $a,b \in \Omega$, $\widetilde{0}_{ab} =  \widetilde{0}_a + W_a \widetilde{0}_b$. 
For  $a \in \Omega$, let $\eta_a := - \widetilde{0}_a$.
Then,$\{W(\eta_a)\}_{a \in \Omega}$ is a gauge cocycle of $\beta$. Note that $\left\lbrace W(\eta_a) \theta_a e(0) = c_{a, 0} e^{\frac{||\eta_a||^2}{2}} e(0)\right\rbrace _{a \in \Omega}$ is a unit of $\beta$. By Proposition \ref{units}, there exists $\chi \in Hom(G,\mathbb{C}^{\times})$ such that \[\chi(a) W(\eta_a) \theta_a e(0) = e(0).\] 
Taking norm in the above equality, we observe that $\chi \in Hom(G,\mathbb{T})$.  Observe that $\{\chi(a)W(\eta_a)\}_{a \in \Omega}$ is a gauge cocycle of $\beta$.  

For $\xi \in Ker(V_{a}^*)$, note that 
\begin{align*}
1 &= \langle e(\xi) | e(0) \rangle \\
&= \langle \chi(a) W(\eta_a) \theta_a e(\xi) | \chi(a) W(\eta_a) \theta_a e(0) \rangle \\
&= \langle \chi(a) c_{a, \xi} e^{ - \frac{|| \eta_a||^2}{2} - \langle \widetilde{\xi}_a | \eta_a \rangle} e(\widetilde{\xi}_a + \eta_a) | e(0) \rangle \\
&=  \chi(a) c_{a, \xi} e^{ - \frac{|| \eta_a||^2}{2} - \langle \widetilde{\xi}_a | \eta_a \rangle} \langle e(\widetilde{\xi}_a + \eta_a) | e(0) \rangle \\
&= \chi(a) c_{a, \xi} e^{ - \frac{|| \eta_a||^2}{2} - \langle \widetilde{\xi}_a | \eta_a \rangle}.
\end{align*}
Thus, $\chi(a) W(\eta_a) \theta_a e(\xi) = e(\widetilde{\xi}_a + \eta_a).$

Suppose $\xi \in Ker(V_a^*) \cap Ker(V_b^*)$ with $a,b \in \Omega$. Choose $x \in \Omega$ such that $a < x$ and $b < x$. This is possible since $a\Omega \cap b\Omega$ is non-empty. Then, $\xi \in Ker(V_a^*) \subset Ker(V_x^*)$, and
\begin{align*}
 e(\eta_x + \widetilde{\xi}_x) &= \chi(x)W(\eta_x) \theta_x e(\xi) \\ &= \chi(a)W(\eta_a)\theta_a e(\xi) \chi(a^{-1}x)W(\eta_{a^{-1}x}) \theta_{a^{-1}x} e(0) \\ &= e(\eta_a + \widetilde{\xi}_a)e(0) =  e(\eta_a + \widetilde{\xi}_a).
\end{align*}
Hence, $\eta_x + \widetilde{\xi}_x = \eta_a + \widetilde{\xi}_a$. Similarly, $\eta_x + \widetilde{\xi}_x = \eta_b + \widetilde{\xi}_b$. Hence, $\eta_a + \widetilde{\xi}_a = \eta_b + \widetilde{\xi}_b$, whenever $\xi \in Ker(V_a^*) \cap Ker(V_b^*)$.

Let  $D=\bigcup_{a \in \Omega} ker ~V_a^*$. Since $V$ is pure, $D$ is dense in $\clh$.   Define 
$ U: D \longrightarrow \mathcal{K}$ by 
\[ U (\xi) := \eta_a + \widetilde{\xi}_a \]
for $\xi \in Ker(V_a^*)$. 
This map is well-defined, as proved earlier. For $a \in \Omega$, the restriction $U|_{Ker(V_a^*)}$ is  continuous.

Fix $\xi, \zeta$ in $D$. Choose $a \in \Omega$ such that $\xi,\zeta \in ker ~V_a^*$. Calculate as follows to observe that
\begin{align*}
e^{ \langle U(\xi) | U(\zeta) \rangle }
&= e^{ \langle \eta_a + \widetilde{\xi}_a | \eta_a + \widetilde{\zeta}_a \rangle} \\
&= \langle e(\eta_a + \widetilde{\xi}_a) | e(\eta_a + \widetilde{\zeta}_a) \rangle\\ 
 &= \langle \chi(a)W(\eta_a)\theta_a e(\xi) | \chi(a)W(\eta_a)\theta_a e(\zeta) \rangle \\
&= \langle e(\xi) | e(\zeta) \rangle \\
&= e^{\langle \xi | \zeta \rangle}. 
\end{align*}
Since $U|_{Ker(V_a^*)}$ is continuous, it follows that there exists an integer $k$ such that for $\xi,\eta \in Ker(V_a^*)$, $ \langle U \xi | U \zeta \rangle = \langle \xi | \zeta \rangle + 2k \pi i$. If $\xi = \zeta$, then $k=0$.  Hence, $U : D \longrightarrow \mathcal{K}$ is an isometry. It is clear that the range of $U$ is dense in $\clk$. Denote the extension of  $U$ again by $U$. Then, $U$ is a unitary. 

We claim that $U$ intertwines $V$ and $W$, i.e. $U V_x = W_x U$ for all $x \in \Omega$. Fix $x \in \Omega$ and $\xi, \zeta \in D$. Choose $a \in \Omega$ such that $x<a$ and $V_x\xi, \zeta \in ker ~V_a^*$. Then $\xi \in Ker(V_{x^{-1}a}^*)$ and we have
\begin{align*}
e^{ \langle UV_x \xi | U \zeta \rangle}
&= e^{ \langle V_x \xi | \zeta \rangle} = \langle e(V_x \xi) | e(\zeta) \rangle \\
&= \langle \chi(a) W(\eta_a) \theta_a e(0 + V_x \xi) | \chi(a) W(\eta_a) \theta_a e(\zeta) \rangle \\
&= \langle \chi(x)W(\eta_x) \theta_x e(0) \chi(x^{-1}a)W(\eta_{x^{-1}a}) \theta_{x^{-1}a}e(\xi) | \chi(a) W(\eta_a) \theta_a e(\zeta) \rangle \\
&= \langle e(0) e(\eta_{x^{-1}a} + \widetilde{\xi}_{x^{-1}a}) | e(\eta_a + \widetilde{\zeta}_a) \rangle \\
&= \langle e(0 + W_x (\eta_{x^{-1}a} + \widetilde{\xi}_{x^{-1}a})) | e(\eta_a + \widetilde{\zeta}_a) \rangle \\
&= \langle e( W_x U \xi) | e (\eta_a + \widetilde{\zeta}_a) \rangle \\
&= e^{ \langle W_x U \xi | U \zeta \rangle}.
\end{align*}  
Hence, there exists an integer $k$ such that for $\xi,\zeta \in D$,  $\langle UV_x \xi | U \zeta \rangle = \langle W_x U \xi | U \zeta \rangle + 2k \pi i$. When $\xi =\zeta=0$, we have $k=0$. Therefore, for $\xi,\eta \in D$ and $x \in \Omega$, we have
$\langle UV_x \xi | U \zeta \rangle = \langle W_x U \xi | U \zeta \rangle$.  Thus, $UV_x = W_xU, $ for all $x \in \Omega$, where $U : \mathcal{H} \longrightarrow \mathcal{K}$ is a unitary. Hence, $V$ and $W$ are unitarily equivalent.

The converse is   omitted as it is straightforward. \hfill $\Box$

Let us discuss the implication of the above theorem to the CCR flows $\alpha^{(A,k)}$. In view of the above theorem, the classification of the CCR flows $\alpha^{(A,k)}$ boils down to the classification of the isometric representation $V^{(A,k)}$. 
We show that, for a general Lie semigroup, the isometric representation $V^{(A,k)}$ remembers $A$, up to a translate, but not necessarily the multiplicity $k$. This is in  contrast to the case of a closed convex cone. 

\begin{rmrk}
Let $A$ be a $P$-space and let $k \in \{1,2,\cdots,\} \cup \{\infty\}$. Suppose $V^{(A,k)}$ is the isometric representation associated to the $P$-space $A$ with multiplicity $k$. A quick explanation, based on groupoids, of the fact that $V^{(A,k)}$ remembers $A$, up to a 
translate,  is described below. 

 The   isometric representations $V^{(A,k)}$, as $A$ and $k$ vary, share the common property that they have commuting range projections. It was proved in \cite{Sundar_Ore} that isometric representations of $P$ with commuting range projections are in $1$-$1$ correspondence with representations of the $C^{*}$-algebra of a universal groupoid which we denote by $\mathcal{G}_{u}$. Moreover, the unit space of $\mathcal{G}_{u}$ is  made of $P$-spaces with an appropriate topology and two elements of the unit space are in the same orbit if and only if one is a translate of other. 
 
 If we appeal to this bijective correspondence, the isometric representation $V^{(A,k)}$ corresponds to the induced representation at the point $A \in \mathcal{G}_{u}$ with multiplicity $k$. It is well known that for groupoids induced representations at  points  not in the same orbit give rise to disjoint representations. Thus, it follows that $V^{(A,k)}$ remembers $A$ up to a translate. 
  To argue that $V^{(A,k)}$ need not remember $k$, it is necessary to pass to a transformation groupoid which is equivalent to $\mathcal{G}_{u}$ and then appeal to Prop. \ref{induced repn of stabiliser} that is proved below. 
 
 However,   the facts that we alluded to in the above paragraphs were not made  explicit neither in \cite{Sundar_Ore} nor in \cite{Anbu_Sundar} where only the case of a cone  was treated. 
  Due to this and also to keep the exposition simpler and the paper fairly self contained,   we choose to directly work with the relevant transformation groupoid, equivalently an ordinary crossed product and avoid any mention of groupoids. We make use of a simple dilation trick. This, in our view is more down-to-earth and makes the paper easier to read.
\end{rmrk}
Let us first collect a few essential things from crossed products that we need. What follows regarding crossed products is  well known and we do not claim originality. We include some details for completeness. 

Let $G$ be a locally compact second countable Hausdorff topological group. Suppose $Y$ is a locally compact second countable Hausdorff space with a continuous left $G$-action. Fix a point $x \in Y$.  Let $k \in \{1,2,\cdots\} \cup \{\infty\}$ be given and let $\clk$ be a Hilbert space of dimension $k$. Consider the Hilbert space $\clh_{k}:=L^{2}(G,\clk)$. Let $\lambda_{G}^{(k)}$ be the left regular representation of $G$ on $\clh_k$. Define a representation $M^{(x,k)}:C_0(Y) \to B(\clh_k)$ by the following formula.
\[
M^{(x,k)}(f)(\xi)(s):=f(s.x)\xi(s)\]
for $\xi \in \clh_k$. 

It is clear that $(M^{(x,k)},\lambda_G^{(k)})$ is a covariant representation of the dynamical system $(C_0(Y),G)$. Observe that the spectral measure of $M^{(x,k)}$ is supported on the orbit containing $x$. Consequently, $M^{(x,k)}$ and $M^{(y,\ell)}$ are disjoint if $x$ and $y$ are in different orbits. The proof of the following Proposition is essentially contained in the commutative diagram given in Lemma 8.26 of \cite{Williams_Dana} and follows by an application of Mackey's imprimitivity theorem. 

\begin{ppsn}
\label{induced repn of stabiliser}
Let $x \in Y$ and let $H$ be the stabiliser of $x$. Let $k,\ell \in \{1,2,\cdots,\} \cup \{\infty\}$ be given. With the foregoing notation, the  following are equivalent.
\begin{enumerate}
\item[(1)] The covariant representations $(M^{(x,k)},\lambda_G^{(k)})$ and $(M^{(x,\ell)},\lambda_G^{(\ell)})$ are unitarily equivalent.
\item[(2)] The left regular representations $\lambda_H^{(k)}$ and $\lambda^{(\ell)}_{H}$ of the subgroup $H$ are unitarily equivalent. 
\end{enumerate}
\end{ppsn}
\textit{Proof.} Assume that $(1)$ holds. Let $B(Y)$ be the algebra of bounded Borel functions on $Y$. Denote the extension of $M^{(x,k)}$ (and $M^{(x,\ell)}$) to $B(Y)$, obtained via the Riesz representation theory, by $M^{(x,k)}$ (and $M^{(x,\ell)}$) itself. 
Since $G/H$ and $Y$ are Polish spaces, it follows that  the map 
\[
G/H \ni sH \to sx \in Y\]
is a Borel embedding. Via this embedding, we can consider a bounded Borel function on $G/H$ as a bounded Borel function on $Y$ by declaring its value outside $G/H$ to be zero. This way, we embedd $C_0(G/H)$ inside $B(Y)$. 

Then, the restiction of   the covariant representations $(M^{(x,k)}, \lambda_G^{(k)})$ and $(M^{(x,\ell)},\lambda_G^{(\ell)})$ to the dynamical system $(C_0(G/H),G)$ are unitarily equivalent. Mackey's imprimitivity theorem states that covariant representations of $(C_0(G/H),G)$ are in bijective correspondence with unitary representations of $H$. But, the covariant representations $(M^{(x,k)},\lambda_G^{(k)})$ and $(M^{(x,\ell)},\lambda_G^{(k)})$ of $(C_0(G/H),G)$ correspond exactly to the left regular representation $\lambda_H^{(k)}$ and $\lambda_H^{(\ell)}$ respectively. Hence $(2)$ holds.
This proves the implication $(1) \implies (2)$. 

 Let $M^{(x,k)}\rtimes \lambda_G^{(k)}$ and $M^{(x,\ell)} \rtimes \lambda_G^{(\ell)}$ be the representations of the crossed product $C_0(G/H) \rtimes G$ that correspond, (i.e. the integrated form),  to the  covariant representations $(M^{(x,k)}, \lambda_G^{(k)})$ and $(M^{(x,\ell)},\lambda_G^{(\ell)})$ respectively. For a representation $\omega$ of $C^{*}(H)$, let $Ind(\omega)$ be the representation of $C_0(G/H) \rtimes G$ obtained via the Rieffel induction using the imprimitivity module that provides the Morita equivalence between the $C^{*}$-algebras $C_0(G/H) \rtimes G$ and $C^{*}(H)$.

Let $M(C_0(G/H) \rtimes G)$ be the multiplier algebra of the crossed product $C_0(G/H) \rtimes G$. 
Let $k \rtimes j_G:C_0(Y) \rtimes G \to M(C_0(G/H) \rtimes G)$ be the homomorphism given in Lemma 8.26 of \cite{Williams_Dana}. Thanks to Lemma 8.26 of \cite{Williams_Dana}, we have
\[M^{(x,k)}\rtimes \lambda_G^{(k)}=Ind(\lambda_H^{(k)})\circ (k \rtimes j_G)~~\textrm{and}~~ M^{(x,\ell)}\rtimes \lambda_G^{(\ell)}=Ind(\lambda_{H}^{(\ell)})\circ (k \rtimes j_G).\]
The  implication $(2) \implies (1)$ follows from the previous equality. \hfill $\Box$

Hereafter, assume that $G$ is a connected Lie group and $P$ is a Lie semigroup with dense interior $\Omega$. Assume that $PP^{-1}=P^{-1}P=G$. Let $\mathcal{C}(G)$ be the set of closed subsets of $G$ equipped with the Fell topology. 
Equipped with the Fell topology, $\mathcal{C}(G)$ is a compact metrisable space. Let us recall the convergence of sequences in $\mathcal{C}(G)$. 

Let $d$ be a metric on $G$ which is compatible with the topology of $G$.  For a closed subset $A$ and $x \in G$, let 
$
d(x,A):=\inf\{d(x,y): y \in A\}.$
For a sequence $(A_n)$ of closed subsets of $G$, define
\begin{align*}
\liminf A_n:&=\{x \in G: \limsup d(x,A_n)=0\}, \\
\limsup A_n:&=\{x \in G: \liminf d(x,A_n)=0\}.
\end{align*}
Then $(A_n)$ converges in $\mathcal{C}(G)$ if and only if $\limsup A_n=\liminf A_n$. In that case,  the sequence $(A_n)$ converges to $\limsup A_n$. Note that $G$ acts continuously on $\mathcal{C}(G)$ by left translations. 

Define 
\begin{align*}
Y_u:&=\{A \in \mathcal{C}(G): A \neq \emptyset, AP^{-1} \subset A\} \\
X_u:&=\{A \in \mathcal{C}(G): e \in A, AP^{-1}\subset A\} \\
X_u^{(0)}:&=\{A \in Y_u: A \cap \Omega \neq \emptyset\}.
\end{align*}
Note that $X_u$ is a compact subset of $\mathcal{C}(G)$. Also, $X_u^{(0)}$ is an open subset of $Y_u$. Observe that $X_u^{(0)} \subset X_u$.  Note that $Y_u$ is invariant under the action of $G$. Moreover, $P X_u \subset X_u$ and
$\Omega X_u \subset X_u^{(0)}$. 

Let $G^{op}:=G$ be the opposite group and consider the preorder on $G^{op}$ induced by the semigroup $P$. 
Choose a sequence $(s_n)$ in $\Omega$ such that $\{s_n:n\geq 1\}$ is cofinal in $G^{op}$. We can also assume that $s_{n+1}>s_n$, i.e. $s_{n+1}s_n^{-1} \in \Omega$. 
We claim that 
\begin{equation}
\label{local compactness}
Y_u=\bigcup_{n=1}^{\infty}s_n^{-1}X_u=\bigcup_{n=1}^{\infty}s_n^{-1}X_u^{(0)}.
\end{equation}
Let $A \in Y_u$ be given. Pick a point $x \in A$. Since $\{s_n: n\geq 1\}$ is cofinal in $G^{op}$, there exists a natural number $n$ such that $s_nx \in \Omega$. This implies that $s_nA \cap \Omega \neq \emptyset$. Hence, $s_nA \in X_u^{(0)}$. Consequently,
$A=s_n^{-1}(s_nA) \in s_n^{-1}X_u^{(0)}$. This proves the claim. 

Since $s_{n+1}s_n^{-1}\in \Omega$, it is clear that $s_{n}^{-1}X_u^{(0)}$ is an increasing sequence of open sets. Moreover, the sets $s_n^{-1}X_u^{(0)}$ have compact closure. Eq. \ref{local compactness} implies that $Y_u$ is locally compact. Since $\mathcal{C}(G)$ is a compact metrisable space, it follows that $Y_u$ is second countable and Hausdorff.

The dynamical system that we make use of is $(C_0(Y_u),G)$. The dynamical system $(C_0(Y_u),G)$ was first considered by Hilgert and Neeb in \cite{Hilgert_Neeb}. View $L^{\infty}(G)$ as the dual of $L^{1}(G)$ and endow $L^{\infty}(G)$ with the weak$^*$-topology. For $A \in Y_u$, let $1_A$ be the indicator function of $A$. Thanks to Prop. II.13 of \cite{Hilgert_Neeb}, the fact that $X_u$ is compact and Eq. \ref{local compactness}, it follows that the map
\[
Y_u \ni A \to 1_A \in L^{\infty}(G)
\]
is a continuous embedding. For $f \in C_c(G)$, let $\widetilde{f}:Y_u \to \mathbb{C}$ be defined by the following equation
\begin{equation}
\label{generating functions}
\widetilde{f}(A):=\int f(x)1_A(x)dx.
\end{equation}

Fix $f \in C_{c}(G)$. We claim that $\widetilde{f} \in C_c(Y_u)$. The continuity of $\widetilde{f}$ follows from the fact that the map $Y_u \ni A \to 1_A \in L^{\infty}(G)$ is continuous. Let $K$ be the support of $f$. Observe that $G$ is the increasing union of open sets $s_n^{-1}\Omega$. Hence $K \subset s_n^{-1}\Omega$ for some $n$. It is clear that 
$\widetilde{f}(A)=0$ if $A \cap s_n^{-1}\Omega=\emptyset$. This means that $\widetilde{f}$ vanishes outside $s_n^{-1}X_u^{(0)}$ which has compact closure. This proves the claim. 

\begin{rmrk}
A straightforward application of the Stone-Weierstrass theorem implies that $\{\widetilde{f}:f \in C_c(G)\}$ generates $C_0(Y_u)$.

\end{rmrk}

Let $V$ be an isometric representation of $P$ on a Hilbert space $\clh$. By the minimal unitary dilation of $V$, we mean a strongly continuous unitary representation $U:=\{U_x\}_{x \in G}$, of $G$,  on a Hilbert space $\clk$ containing $\clh$ as a closed subspace such that 
\begin{enumerate}
\item[(1)] for $a \in P$ and $\xi \in \clh$, $U_a\xi=V_a\xi$, and
\item[(2)] the union $\bigcup_{a \in P}U_a^{*}\clh$ is dense in $\clk$.
\end{enumerate}
It is clear that the minimal unitary dilation is unique up to unitary equivalence. For its existence,  we refer the reader to \cite{Laca95}.

Let $V^{(1)}$ and $V^{(2)}$ be isometric representations of $P$ on Hilbert spaces $\clh_1$ and $\clh_2$ respectively. For $i=1,2$, let $U^{(i)}$ be the minimal unitary dilation of $V^{(i)}$ and suppose that $U^{(i)}$ acts on $\clk_i$. 
For $x \in G$, let $E^{(i)}_{x}$ be the orthogonal projection onto $U_x^{(i)}\clh_i$. The proof of the following proposition is quite elementary and hence omitted. 

\begin{ppsn}
\label{dilation trick}
Keep the foregoing notation. The following are equivalent.
\begin{enumerate}
\item[(1)] The isometric representations $V^{(1)}$ and $V^{(2)}$ are unitarily equivalent.
\item[(2)] There exists a unitary $U:\clk_1 \to \clk_2$ such that for every $x \in G$, 
\[
UU^{(1)}_x U^{*}=U^{(2)}_x ~~\textrm{and~~}UE^{(1)}_xU^{*}=E^{(2)}_{x}.\]
\end{enumerate}
\end{ppsn}

Let $A_1$ and $A_2$ be $P$-spaces. Suppose $\clk_1$ and $\clk_2$ are Hilbert spaces of dimension $k_1$ and $k_2$ respectively. 
For $i=1,2$, set $\clh_i:=L^{2}(A_i,\clk_i)$. Denote the isometric representation of $P$ on $\clh_i$ associated to the $P$-space $A_i$
with multiplicity $k_i$ by $V^{(i)}$. Fix $i \in \{1,2\}$. Let $\lambda^{(k_i)}$ be the left regular representation of $G$ on $L^{2}(G,\clk_i)$. 
We view $\clh_i:=L^{2}(A_i,\clk_i)$ as a closed subspace of $L^{2}(G,\clk_i)$ in the obvious way. 

It is clear the for $x \in P$, $V^{(i)}_x$ is the compression of $\lambda^{(k_i)}_{x}$ onto $\clh_i$. By replacing $A_i$ with a right translate of $A_i$, we can assume without loss of generality
that $e \in A_i$. Then $P \subset A_i$. Since $G=\Omega^{-1}\Omega$, it follows that $G=\bigcup_{a \in \Omega}a^{-1}A_i$. This has the consequence that the union
$\displaystyle \bigcup_{a \in \Omega}\lambda_{a^{-1}}^{(k_i)}\clh_i$ is dense in $L^{2}(G,\clk_i)$. Hence, $\lambda^{(k_i)}$ is the minimal unitary dilation of $V^{(i)}$. 

Let $M$ be the multiplication representation of $L^{\infty}(G)$ on $L^{2}(G,\clk_i)$. For $x \in G$, let $E_{x}^{(i)}$ be the orthogonal projection onto $\lambda^{(k_i)}_x\clh_i$. Then, it is clear 
that $E^{(i)}_x=M(1_{xA_i})$. 
For $i \in \{1,2\}$, let $B_i:=A_{i}^{-1}$. Recall the covariant representations $(M^{(B_i,k_i)},\lambda_{G}^{(k_i)})$ of the dynamical system $(C_0(Y_u),G)$ explained before Prop. \ref{induced repn of stabiliser}. With the foregoing notation, we have the following.

\begin{ppsn}
\label{algebraise}
The following are equivalent.
\begin{enumerate}
\item[(1)] The isometric representations $V^{(1)}$ and $V^{(2)}$ are unitarily equivalent. 
\item[(2)] The covariant representations $(M^{(B_1,k_1)},\lambda_G^{(k_1)})$ and $(M^{(B_2,k_2)},\lambda_G^{(k_2)})$ are unitarily equivalent. 
\end{enumerate}
\end{ppsn}
\textit{Proof.} In what follows, we simply denote $\lambda_G^{(k_i)}$ by $\lambda^{(k_i)}$. Fix $i \in \{1,2\}$. Let $f \in C_{c}(G)$ be given. First, observe that for $\xi,\eta \in L^{2}(G,\clk_i)$, we have

\begin{align*}
\Big \langle \big(\int f(x)E_x^{(i)}dx\big)\xi|\eta \Big \rangle &= \int f(x) \Big(\int 1_{xA_i}(s)\langle \xi(s)|\eta(s) \rangle ds \Big)dx \\
&= \int \big(\int f(x)1_{sB_i}(x)dx\big) \langle \xi(s)|\eta(s)\rangle ds\\
&=\int \langle \widetilde{f}(sB_i)\xi(s)|\eta(s)\rangle ds.
\end{align*} 
Here, $\widetilde{f}$ stands for the function defined by Eq. \ref{generating functions}. Consequently, we have 
\begin{equation}
\label{Master equation}
\int f(x)E_x^{(i)}dx=M^{(B_i,k_i)}(\widetilde{f})
\end{equation}
for $f \in C_c(G)$. 

Suppose that $(1)$ holds.  Thanks to Prop. \ref{dilation trick} and the fact that the minimal unitary dilation of $V^{(i)}$ is $\lambda^{(k_i)}$, we have a unitary $U:L^{2}(G,\clk_1) \to L^{2}(G,\clk_2)$ such that 
\[
U\lambda_x^{(k_1)}U^{*}=\lambda_x^{(k_2)} ~~\textrm{and~~} UE_{x}^{(1)}U^{*}=E_{x}^{(2)}\]
for every $x \in G$. Eq. \ref{Master equation} implies that for $f \in C_{c}(G)$, $UM^{(B_1,k_1)}(\widetilde{f})U^{*}=M^{(B_2,k_2)}(\widetilde{f})$. But $\{\widetilde{f}: f \in C_c(G)\}$ generates $C_0(Y_u)$ (see Remark \ref{generating functions}). Hence, $U$ intertwines $M^{(B_1,k_1)}$ and $M^{(B_2,k_2)}$. The unitary $U$ already intertwines $\lambda^{(k_1)}$ and $\lambda^{(k_2)}$. Therefore, $(2)$ holds. This proves the implication $(1) \implies (2)$. 

Suppose that there exists a unitary $U:L^{2}(G,\clk_1) \to L^{2}(G,\clk_2)$ that intertwines $(M^{(B_1,k_1)},\lambda_G^{(k_1)})$ and $(M^{(B_2,k_2)},\lambda_G^{(k_2)})$. Appealing to Eq. \ref{Master equation}, we see that 
for $f \in C_c(G)$,
\[
\int f(x)UE_x^{(1)}U^{*}dx=\int f(x)E_x^{(2)}dx.\]
As the above equality holds for every continuous compactly supported function and the maps $G \ni x \to E_x^{(i)} \in B(L^{2}(G,\clk_i))$ are weakly continuous, it follows that for every $x \in G$, $UE_x^{(1)}U^{*}=E_x^{(2)}$. By Prop. \ref{dilation trick}, it follows that $V^{(1)}$ and $V^{(2)}$ are unitarily equivalent. This proves the implication $(2) \implies (1)$. \hfill $\Box$

The following two corollaries are immediate from Theorem \ref{injectivity}, Prop. \ref{algebraise}, Prop. \ref{induced repn of stabiliser} and the discussion preceeding Prop. \ref{induced repn of stabiliser}. 
\begin{crlre}
Keep the foregoing notation. Suppose that $V^{(1)}$ and $V^{(2)}$ are unitarily equivalent. Then there exists $z \in G$ such that $A_1z=A_2$. 
\end{crlre}

Let $A$ be a $P$-space and $k,\ell \in \{1,2,\cdots,\} \cup \{\infty\}$. Let $V^{(A,k)}$ and $V^{(A,\ell)}$ be the isometric representations associated to the $P$-space $A$ with multiplicity $k$ and $\ell$. Denote the 
associated CCR flows by $\alpha^{(A,k)}$ and $\alpha^{(A,\ell)}$ respectively. Let 
\[
G_A:=\{z \in G: Az=A\}.\]
Observe that $G_A$ is a closed subgroup of $G$. 

\begin{crlre}
With the foregoing notation, the  following are equivalent. 
\begin{enumerate}
\item[(1)] The CCR flows $\alpha^{(A,k)}$ and $\alpha^{(A,\ell)}$ are cocycle conjugate.
\item[(2)] The isometric representations $V^{(A,k)}$ and $V^{(A,\ell)}$ are unitarily equivalent.
\item[(3)] The left regular representation of $G_A$ with multiplicity $k$ is unitarily equivalent to the left regular representation of $G_A$ with multiplicity $\ell$.
\end{enumerate}
\end{crlre}

\begin{rmrk}
\begin{enumerate}
\item[(1)] If the semigroup $P$ is abelian, then the CCR flows $\alpha^{(A,k)}$ remember the multiplicity $k$. 
\item[(2)] Let $G:=H_{5}$ be the Heisenberg group of dimension $5$ and let $P:=H_{5}^{+}$ be the Heisenberg subsemigroup consisting of non-negative entries. 
Let \[A:=\Big\{\begin{bmatrix}
                                                                                                         1& x_1 & x_2 & z \\
                                                                                                         0 & 1 & 0 & y_1 \\
                                                                                                         0 & 0 & 1 & y_2 \\
                                                                                                         0 & 0 & 0 & 1
                                                                                                           
                                                                                                                                                                                                                                             \end{bmatrix} \in G: x_1 \geq 0, y_1 \geq 0 \Big\}.\]
  Clearly,  $A$ is a $P$-space and $G_A$ is isomorphic to the three dimensional Heisenberg group $H_3$.

  The Plancherel theorem for $H_3$ states that the left regular representation with multiplicity $1$ disintegrates
  into irreducible representations where one dimensional representations does not occur and each irreducible infinite dimensional representation occur with infinite multiplicity. Consequently, the left regular representation of any multiplicity
  disintegrates the same way as the left regular representation with multiplicity $1$.    Hence, the left regular representation of any multiplicity is unitarily equivalent to the left regular representation with multiplicity $1$. 
  
  This has the implication 
  that the CCR flows $\alpha^{(A,k)}$, as $k$ varies, belong to the same cocycle conjugacy class.

\end{enumerate}
\end{rmrk}

\section{Uncountably many type I examples}

In this section, we produce the promised uncountably many type I CCR flows with index one. For the rest of this paper, the letter $P$ stands for a closed convex cone
in $\mathbb{R}^{d}$, which we assume is spanning and pointed. We also assume $d \geq 2$. 

Suppose $\alpha$ is an $E_0$-semigroup over $P$ on $B(\clh)$. Let $E:=\{E(x)\}_{x \in P}$ be the product system of $\alpha$. Suppose that, for $x \in P$,  $F(x)$ is a non-zero closed subspace of $E(x)$. 
Set $F:=\{F(x)\}_{x \in P}$. We say that $F$ is a \emph{subsystem} of $E$ if for $x,y \in P$, 
\[
F(x+y)=\overline{\textrm{span}\{ST: S \in F(x), T \in F(y)\}}.\]
In other words, $F$ is a product system on its own right.  
Let $\alpha$ be a  spatial $E_0$-semigroup. Denote the set of units of $\alpha$ by $\mathcal{U}_\alpha$. We say that a subsystem $F=\{F(x)\}_{x \in P}$ of $E$ contains $\mathcal{U}_\alpha$  if
for every $u=\{u_x\}_{x \in P} \in \mathcal{U}_{\alpha}$, $u_x \in F(x)$ for $x \in P$. 

\begin{dfn}
Let $\alpha:=\{\alpha_x\}_{x \in P}$ be an $E_0$-semigroup and let $E:=\{E(x)\}_{x \in P}$ be the product system of $\alpha$. Suppose that $\alpha$ is spatial. Let $F:=\{F(x)\}_{x \in P}$ be a subsystem of $E$. 
We say $F$ is the type I part of $E$ or $\alpha$ if it satisfies the following. 
\begin{enumerate}
\item[(1)] The subsystem $F$ contains $\mathcal{U}_\alpha$.
\item[(2)] If $G:=\{G(x)\}_{x \in P}$ is a subsystem of $E$ that contains $\mathcal{U}_\alpha$ then $F(x) \subset G(x)$ for every $x \in P$. 
\end{enumerate}
We say $\alpha$ is type I if the type I part of $E$ is $E$. 
\end{dfn}

\begin{rmrk}
It is not clear whether, in the higher dimensional case, the type I part of a product system always exists. For a spatial $E_0$-semigroup $\alpha$ with the associated  product system $E:=\{E(x)\}_{x \in P}$, we could set 
for $x \in P$
\[
F(x):=\overline{span\{u^{(1)}_{x_1}u^{(2)}_{x_2}\cdots u^{(n)}_{x_n}: u^{(i)} \in \mathcal{U}_{\alpha}, \sum_{i=1}^{n}x_i=x, n \in \mathbb{N}\}}.\]
When $P=[0,\infty)$, $F:=\{F(x)\}_{x \in P}$ is a subsystem of $E$ and it is the type I part of $E$. 
However, in the higher dimensional case, it is unclear whether $F:=\{F(x)\}_{x \in P}$ is a subsystem. In particular, it is not clear whether $F(x+y)\subset F(x)F(y)$. 
For the order induced by the cone is only a partial order and not a total order. 

In this context, we pose the following two questions. 
\begin{enumerate}
\item[(1)] Does the type I part of a product system  exist ?
\item[(2)] Is the field of Hilbert spaces $F:=\{F(x)\}_{x \in P}$, defined above,  a subsystem ?
\end{enumerate}
\end{rmrk}  

We next show that the type I part of a CCR flow exists.   Fix a strongly continuous isometric representation $V$ of $P$ on $\clh$. We assume that $V$ is pure. Let $\alpha^{V}$ be the CCR flow associated to $V$ and let $E$ be the product system of $\alpha^{V}$. 
As before, for $x \in P$, we  identify $E(x)$ with $\Gamma(Ker(V_{x}^{*}))$. 

A closed subspace $\clh_0$ of $\clh$ is said to be invariant under $V$ if $\clh_0$ is invariant under $\{V_x,V_{x}^{*}:x \in P\}$. Let $\clh_0$ be the smallest closed subspace of $\clh$ invariant under $V$ which contains the set $\{\xi_{x}: x \in P, \xi \in \mathcal{A}(V)\}$. Denote the orthogonal projection onto $\clh_0$ by $Q$. Define for $x \in P$, 
\[
F(x):=\overline{span\{e(\eta): Q\eta=\eta, \eta \in Ker(V_x^{*})\}}.\]
Set $F:=\{F(x)\}_{x \in P}$. Then $F$ is a subsystem of $E$. Note that $F$ is isomorphic to the product system of the CCR flow associated to the isometric representation $V$ restricted to $\clh_0$. Recall from \cite{Anbu}, or from Proposition \ref{units}, that for $\xi \in \mathcal{A}(V)$, $\{e(\xi_x)\}_{x \in P}$ is a unit of $E$ and every unit, up to a character, arises this way. Hence $F$ contains the units of $E$. 

\begin{ppsn}
\label{smallest subsystem}
With the above notation, $F$ is the type I part of $E$. 
\end{ppsn}
\textit{Proof.} Let $G:=\{G(x)\}_{x \in P}$ be a subsystem of $E$ which contains the units of $E$. For $x \in P$, let $\theta_{x}:E(x) \to E(x)$ be the orthogonal projection onto $G(x)$. The fact that $G$ is a subsystem implies that  for $u \in E(x)$ and $v \in E(y)$, $\theta_{x+y}(uv)=\theta_{x}(u)\theta_x(v)$. Consequently, $\theta$ is a local projective cocycle. 

By Prop. 6.12 of \cite{injectivity}, it follows that there exists $\xi \in \mathcal{A}(V)$ and a projection $R$ in the commutant of $\{V_x,V_{x}^{*}: x \in P\}$ such that for $x \in P$, $\eta \in Ker(V_x^{*})$, 
\begin{equation}
\label{local cocycle}
\theta_{x}(e(\eta))=e^{-\langle \lambda|x \rangle}e^{\langle \eta|\xi_x \rangle}e(R\eta+\xi_{x}).
\end{equation}
where $\lambda$ is such that $\langle \lambda|x \rangle=\langle \xi_x|\xi_x\rangle$ for $x \in P$. 

Note that $\{e(0)\}_{x \in P}$ is a unit of $E$. Hence $\theta_{x}(e(0))=e(0)$. Eq. \ref{local cocycle} implies that $\xi_{x}=0$. Let $\eta \in \mathcal{A}(V)$ be given. Since $G$ contains the units of $E$,  we have $\theta_{x}(e(\eta_x))=e(\eta_x)$.  Eq. \ref{local cocycle} implies that $R\eta_x=\eta_x$. Consequently, $Q \leq R$. Now it is clear that for $x \in P$, $F(x) \subset G(x)$. This completes the proof. \hfill $\Box$

\begin{rmrk}
\label{type one vs irred}
For a pure isometric representation $V$ on $\clh$,  the CCR flow $\alpha^{V}$  is type I if and only if the smallest  closed  subspace invariant under $V$ that contains  the set $\{\xi_x: x \in P, \xi \in \mathcal{A}(V)\}$ is $\clh$. 

In particular, if the isometric representation has a non-zero additive cocycle and is irreducible, in the sense that the commutant of $\{V_x,V_x^{*}:x \in P\}$ is trivial, then $\alpha^{V}$ is type I. 
\end{rmrk}

Now, we construct the desired type I examples. Let $P^{*}$ be the dual cone of $P$, i.e. 
\[P^*:= \left\lbrace y \in \mathbb{R}^d :  \langle x| y \rangle \geq 0, \forall x \in P \right\rbrace.\]
Then, $P^{*}$ is spanning and pointed. Choose $e \in Int(P^{*})$ of norm one. The vector $e$ will be fixed for the rest of this paper. 

Let $N \subset (span\{e\})^{\perp}$ be a discrete subgroup of rank $d-1$. We claim that for a closed subset $F \subset P$, $F+N$ is closed in $\mathbb{R}^d$. If $x \in \overline{F+N}$,  then there exists a sequence $(p_n + m_n) \in F+N$, with $p_n \in F$ and $m_n \in N$,  converging to $x$. Then, \[\langle p_n|e \rangle=\langle p_n+m_n | e \rangle \to \langle x|e \rangle.\] Thus, there exists $c>0$ such that $\langle p_n|e \rangle \leq c$ for every $n$. 
Let $S:= \left\lbrace y \in F : \langle y | e \rangle \leq c \right\rbrace$. 
By Lemma I.1.6 of \cite{Faraut}, $S$ is compact. Since $(p_n)$ is in $S$,  there exists $p \in S$  and a subsequence $(p_{n_k})$ such that $p_{n_k}$ converges to $p$. Then, $m_{n_k}$ converges to $m=x-p \in N$. Thus, $p_{n_k} + m_{n_k} \longrightarrow p+m = x$ and $x \in F+N$. Hence, $F+N$ is closed

Let $\phi : \mathbb{R}^d \ni x \to x+N \in \mathbb{R}^d/N$ be the quotient map. Let 
\begin{align*}
G^{N}:&=\phi (\mathbb{R}^d) = \mathbb{R}^d/N = \mathbb{R} \times \mathbb{T}^{d-1}\\
Q^{N}:&=\phi(P).
\end{align*}
Since $F+N$ is closed whenever $F$ is a closed subset of $P$, it follows that $Q^{N}$ is closed in $G^{N}$. Moreover, the map $\phi: P \to Q^{N}$ is a closed and hence a quotient map. Also, since $P \subset L(Q^N)$, the Lie-wedge of $Q^N$, we have $Q^N = \phi(P) \subset \phi(L(Q^N)) \subset Q^N$, and hence $Q^N$ is a Lie semigroup of $G^N$.

\begin{lmma}
\label{irred}
Keep the foregoing notation. The intersection $Q^{N} \cap -Q^{N}=\{0\}$. 
\end{lmma}
\textit{Proof.} It suffices to show that $P \cap N=\{0\}$. Note that 
$P \cap N \subset P \cap (span\{e\})^{\perp}$. But $e \in Int(P^{*})$. By Prop. I.1.4 of \cite{Faraut}, we have $P \cap (span\{e\})^{\perp}=\{0\}$. Hence, $P \cap N=\{0\}$. This completes the proof. \hfill $\Box$

 Define $W^N: Q^N \longrightarrow L^2(Q^N)$ by 
\[ W^N_y f(z) := \begin{cases}
f(z - y) \quad \textrm{if } z - y \in Q^N, \\
0 \quad \textrm{otherwise}.
\end{cases}\]

 Define $V^N: P \longrightarrow L^2(Q^N)$ by
\[ V^N_x := W^N_{\phi(x)} \]
Note that $V^N$ is an isometric representation of $P$ on $L^2(Q^N)$. We denote the CCR flow associated to $V^{N}$ by $\alpha^{N}$. 
We will use the notation $V^N$, $W^N$ and $\alpha^{N}$ throughout the rest of this section. 

\begin{ppsn}
\label{index one example}
Let $N \subset (span\{e\})^{\perp}$ be a discrete subgroup of rank $d-1$. Then, we have the following.
\begin{enumerate}
\item[(1)] 	The CCR flow  $\alpha^N$ has index one.
\item[(2)] The CCR flow $\alpha^{N}$ is type I.
\end{enumerate} 
\end{ppsn}

\textit{Proof:} If $\xi$ is an additive cocycle of $W^N$, then it is easy to check that $\left\lbrace \eta_x = \xi_{\phi(x)} \right\rbrace _{x \in P}$ is an additive cocycle of $V^N$. Conversely, suppose $\eta = \left\lbrace \eta_x \right\rbrace _{x \in P}$ is an additive cocycle of $V^N$.  Define $\xi: Q^{N} \to L^{2}(Q^{N})$ as follows. For $\widetilde{y}:=y+N \in Q^{N}$, with $y \in P$, set \[\xi_{\widetilde{y}}:=\eta_y.\]

We claim that $\xi$ is well defined. Suppose $y+N=x+N$, with $x,y \in P$, then $V^N_y = W^N_{\widetilde{y}} = W^N_{\widetilde{x}} = V^N_x$ and $\eta_{x+y} = \eta_x + V^N_x \eta_y = \eta_y + V^N_x \eta_x$.
 Hence, \[\eta_x - \eta_y = V^N_x(\eta_x - \eta_y)\] which implies $0= V_x^{N*}(\eta_x - \eta_y) = \eta_x - \eta_y$. Thus, if $\widetilde{x} = \widetilde{y}$, then $ \xi_ {\widetilde{x}} = \eta_x = \eta_y = \xi_{\widetilde{y}}$ and $\xi$ is well-defined.  
  This proves the claim.     Since $\eta$ is continuous and descends onto the quotient space,  it follows that $\xi$ is continuous. It is straightforward to check that $\xi$ is an additive cocycle. 
    
    Thus, every additive cocycle of $V^{N}$ is a pullback of an additive cocycle of $W^{N}$. 
  By Prop. \ref{additive cocycle}, Prop. \ref{sec 3 result} and Theorem \ref{when is boundary}, we have $\mathcal{A}(W^N) =\left\lbrace \left\lbrace  \lambda 1_{(Q^N \backslash \widetilde{x} Q^N)} \right\rbrace_{ \widetilde{x} \in Q} | \lambda \in \mathbb{C} \right\rbrace $. Thus, $\mathcal{A}(V^N) = \left\lbrace \left\lbrace \lambda 1_{( Q^N \backslash \widetilde{x}Q^N)} \right\rbrace_{x \in P} | \lambda \in \mathbb{C} \right\rbrace $. Consequently, the CCR flow $\alpha^{N}$ has index one, i.e. $Ind(\alpha^N) = 1$.

Thanks to Remark \ref{type one vs irred}, to show that $\alpha^{N}$ is type I, it suffices to show that $V^{N}$ is irreducible. But $V^{N}$ is the pullback of $W^{N}$ and consequently, it suffices to show that $W^{N}$ is irreducible. 
Let \[\mathcal{W}:=C^*\Big \{\int f(y)W_{y}^{N}dy: f \in L^{1}(Q^N)\Big \}.\] The $C^{*}$-algebra $\mathcal{W}$ is called the Wiener-Hopf algebra associated to the Lie semigroup $Q^{N}$. It is clear that 
the commutant of $\{W_y^{N}, W_y^{N*}: y \in Q^{N}\}$ coincides with the commutant of $\mathcal{W}$. By Lemma \ref{irred} and by Theorem IV.11 of \cite{Hilgert_Neeb}, it follows that $\mathcal{W}$ contains the algebra of compact operators on $L^{2}(Q^N)$. 
Hence, the commutant of $\mathcal{W}$ is trivial.  This proves that $W^{N}$ is irreducibe and hence the proof. \hfill $\Box$

For a closed subgroup $H$ of $\mathbb{R}^{d}$, set 
\[
H^{\perp}:=\{\xi \in \mathbb{R}^d: e^{i \langle x|\xi \rangle}=1, \forall x \in H\}.
\]

\begin{ppsn} \label{N_1 = N_2}
	Let $N_1$ and $N_2$ be discrete subgroups of $(span\{e\})^{\perp}$ of rank $d-1$.  Let $\phi_i : \mathbb{R}^d \longrightarrow \mathbb{R}^d / N_i$ be the quotient map for $i =1,2$.  Then, $V^{N_1}$ is unitarily equivalent to $V^{N_2}$ iff $N_1 = N_2$.
\end{ppsn}
\textit{Proof:} For $i=1,2$, let $U^{(i)}$ be the minimal unitary dilation of $V^{N_i}$ and let $\lambda^{(i)}$ be the left regular representation of $\mathbb{R}^{d}/N_i$ on $L^{2}(\mathbb{R}^{d}/N_i)$.  Let us assume that $V^{N_1}$ is unitarily equivalent to $V^{N_2}$. Then, the minimal unitary dilations $U^{(1)}$ and $U^{(2)}$ are unitarily equivalent. However, $V^{N_i}$ is the pullback of $W^{N_i}$ and, as observed in Section 5, the minimal unitary dilation of $W^{N_i}$ is $\lambda^{(i)}$. Therefore, for $x \in \mathbb{R}^{d}$, 
we have $U^{(i)}_x=\lambda^{(i)}_{\phi_i(x)}$. Since $U^{(1)}$ and $U^{(2)}$ are unitarily equivalent, for every $x \in \mathbb{R}^{d}$, $U^{(1)}_x$ and $U^{(2)}_x$ have the same spectrum. 

Equating spectrums, we get for $x \in \mathbb{R}^{d}$, 
\[ \sigma(  \lambda^{(1)}_{\phi_1(x)})  =\overline{ \left\lbrace e^{ i \langle x| \xi \rangle} | \xi \in N_1^\bot \right\rbrace} =  \sigma (\lambda^{(2)}_{\phi_2(x)})  = \overline{\left\lbrace e^{ i \langle x| \eta \rangle} | \eta \in N_2^\bot \right\rbrace}. \]
The above equality forces that if $x \in N_1$,  for $\eta \in N_2^{\perp}$, $ e^{i \langle x | \eta \rangle }=1$. In otherwords, we have $N_1 \subset (N_2^{\perp})^{\perp}=N_2$. Similarly, $N_2 \subset N_1$ and hence $N_1=N_2$.  This completes the proof. \hfill $\Box$.

\begin{ppsn} \label{not pullback}
	Let $N$ be a discrete subgroup of $(span\{e\})^{\perp}$ of rank $d-1$. The  isometric representation $V^N$ is not unitarily equivalent to a pullback of any one-parameter isometric representation. 
	\end{ppsn}
\textit{Proof:} Let $\{S_t\}_{t \geq 0}$ be the shift semigroup on $L^{2}(0,\infty)$. Suppose the representation $V^{N}$ is unitarily equivalent to a pullback of a one parameter isometric representation. Since $V^{N}$ is irreducible, it follows that there exists a homomorphism $\psi : P \longrightarrow [0, \infty)$ and a unitary operator $X:L^{2}(Q^N) \to L^{2}(0,\infty)$ such that  $XV^N_xX^* = S_{\psi(x)}$. Since $\psi:\mathbb{R}^{d} \longrightarrow \mathbb{R}$ is a homomorphism, there exists $\mu \in \mathbb{R}^d$ such that $\psi(x) = \langle x| \mu \rangle$. Clearly $\mu \neq 0$. 

Let $U:=\{U_x\}_{x \in \mathbb{R}^{d}}$ be the minimal unitary dilation of $V^{N}$, which we observed in Prop. \ref{N_1 = N_2}, is the pullback of the left regular representation, $\lambda^{N}$, of $\mathbb{R}^{d}/N$ via the quotient map $\phi: \mathbb{R}^{d} \ni x \to x+N \in \mathbb{R}^{d}/N$. On the other hand, the minimal unitary dilation of $\{S_t\}_{t \geq 0}$ is the left regular representation, $\{\lambda_t\}_{t \in \mathbb{R}}$ of $\mathbb{R}$ on $L^{2}(\mathbb{R})$. 

Since, $V^{N}$ is unitarily equivalent to $\{S_{\langle x|\mu \rangle}\}_{x \in P}$, their minimal unitary dilations are equivalent. This implies that for every $x \in \mathbb{R}^{d}$, $U_x:=\lambda^{N}_{\phi(x)}$ and $\lambda_{\langle x|\mu \rangle}$ have the same spectrum. 
Equating their spectrums, we get for every $x \in \mathbb{R}^{d}$,  
\[ \overline{ \left\lbrace e^{i \langle x| \xi \rangle} : \xi \in N^\bot \right\rbrace }= \left\lbrace e^{i t   \langle x| \mu \rangle } : t \in \mathbb{R} \right\rbrace.
\]
The above equality implies that $(span\{\mu\})^{\perp} \subset (N^{\perp})^{\perp}=N$ which is a contradiction. Hence the proof. \hfill $\Box$

\begin{rmrk}
Combining Theorem \ref{injectivity}, Prop. \ref{index one example}, Prop. \ref{N_1 = N_2} and Prop. \ref{not pullback}, we see that there are uncountably many type I CCR flows with index one which are not pullbacks of one parameter CCR flows.
By considering the representation $V^{N} \otimes 1$ on $L^{2}(Q^N) \otimes \clk$, where $\clk$ is a separable Hilbert space of dimension $k$, we find that there are uncountably many type I CCR flows with index $k$ for any $k \in \{1,2,\cdots,\} \cup \{\infty\}$. 
Unlike the case $k=1$, for $k \geq 2$, the representation $V^{N} \otimes 1$ is not irreducible. 

This raises the following  question. 
Does there exists an isometric representation $V$ such that $V$ is irreducible but the associated CCR flow $\alpha^{V}$ has index strictly greater than one ?
\end{rmrk}

\bibliography{references}
 \bibliographystyle{amsplain}

 \nocite{Sundar_Notes}
 \nocite{Vasanth}

 \vspace{2.5mm}
 
\noindent
{\sc Piyasa Sarkar}
(\texttt{piyasa10@gmail.com})\\
{\footnotesize Institute of Mathematical Sciences (HBNI), CIT Campus, \\
Taramani, Chennai, 600113, Tamilnadu, India.}\\

\noindent
{\sc S. Sundar}
(\texttt{sundarsobers@gmail.com})\\
         {\footnotesize  Institute of Mathematical Sciences (HBNI), CIT Campus, \\
Taramani, Chennai, 600113, Tamilnadu, INDIA.}\\

\end{document}